\magnification 1200
\input amstex
\documentstyle{amsppt}
\vsize=7.85in
\topmatter

\rightheadtext{Rational dependence of CR
mappings}
\leftheadtext{M. S. Baouendi,  P.
Ebenfelt, and L. P. Rothschild}
\title Rational dependence of
smooth and analytic CR mappings on
their jets\endtitle
\author M. S. Baouendi,  P. Ebenfelt, and Linda
Preiss Rothschild
\endauthor
\address Department of Mathematics, 0112, University of
California at San Diego, La Jolla, CA
92093-0112\endaddress
\email sbaouendi\@ucsd.edu,
lrothschild\@ucsd.edu\endemail
\address Department of Mathematics, Royal
Institute of Technology, 100 44 Stockholm,
Sweden\endaddress
\email ebenfelt\@math.kth.se\endemail
\thanks The first and the third authors are
partially supported by National Science
Foundation grant DMS 98-01258. The second author
is partially supported by a grant from the Swedish Natural
Science Research Council.\endthanks
\loadeufm

\def\Rk{\text{\rm Rk}}

\def\det{{\text{\rm det}}}

\def\a {\alpha}


\define \codim{\text{\rm codim }}

\define \bR{\Bbb R}
\define \bC{\Bbb C}



\define \Aut{\text{\rm Aut}}

\def\rk{\text{\rm rk}}

\define\Span{\text{\rm span\,}}
\endtopmatter

\document

\heading \S 0. Introduction\endheading

Let $M\subset\bC^N$ and $M'\subset\bC^{N'}$
be two smooth ($C^\infty$) generic
submanifolds with $p_0\in M$ and $p'_0\in
M'$. We shall consider holomorphic mappings
$H\:(\bC^N,p_0)\to (\bC^{N'},p'_0)$, defined
in a neighborhood of $p_0\in\bC^N$, such
that
$H(M)\subset M'$ (and,
more generally, smooth CR mappings
$(M,p_0)\to (M',p'_0)$; see below). We shall always
work under the assumption that $M$ is of
finite type at
$p_0$ in the sense of Kohn and
Bloom--Graham, and that 
$M'$ is finitely nondegenerate at $p_0'$ 
(see \S1 for precise definitions). More
precisely, we shall assume that $M'$ is
$\ell_0$-nondegenerate at $p_0$, for some
integer $\ell_0\geq 0$. (Recall that for a
real hypersurface, $1$-nondegeneracy at a point is equivalent to Levi 
nondegeneracy at  
that point. Also, as is further explained in \S2.2, the notion of finite
nondegeneracy for real-analytic generic submanifolds is 
intimately related to the notion of holomorphic nondegeneracy defined
below in this section.) 

We denote by $J$ the complex structure map
on $T\bC^N$. Recall that for $p\in M$,
$T^c_{p}M$ denotes the {\it complex tangent
space} to $M$ at
$p$, i.e. the largest $J$-invariant
subspace of $T_pM$, the tangent space of
$M$ at $p$. A smooth mapping $H\:M\to
M'$ is called {\it CR} if its tangent map
$dH$ maps $T^c_pM$ into $T^c_{H(p)}
M'$ for every $p\in M$. A CR mapping $H\:
M\to M'$ is called {\it CR submersive} at
$p$ if
$dH$ maps $T^c_pM$ onto $T^c_{H(p)}M'$. A
holomorphic mapping $H$ sending $M$ into
$M'$ is called CR submersive if its
restriction to $M$ is. To a smooth
CR mapping $H\:M\to M'$, and $p_0\in M$, one
may associate a unique formal (holomorphic) 
power series
mapping 
$$
\hat H(Z)\sim\sum_\alpha a_\alpha
(Z-p_0)^\alpha,\quad
a_\alpha\in\bC^{N'},\tag0.1
$$
which sends $(M,p_0)$ into $(M',p_0')$.
(See \S1 for details, definitions, and discussion.) 
If $H$ extends holomorphically to a
neighborhood of $p_0$ in $\bC^N$, then
$\hat H(Z)$ is the Taylor series of $H$ at
$p_0$.  For
$p\in\bC^N$ and $p'\in\bC^{N'}$, we shall
denote by
$J^k(\bC^N,\bC^{N'})_{(p,p')}$ the {\it jet
space} of order $k$ of holomorphic 
mappings $(\bC^N,p)\to (\bC^{N'},p')$. (See
\S2 for further details.) 

Our main
result (which will also be given in a
slightly more precise and general form in
Theorem 2.1.5) gives rational dependence of
the formal power series mapping associated to a CR
submersive mapping, and more generally of a  
formal CR submersive power series mapping 
$(M,p_0)\to 
(M',p_0')$ (see \S1 for precise definitions), 
on its jet of a
predetermined order. 

\proclaim{Theorem 1} Let $M$ and $M'$
be smooth generic submanifolds through
$p_0\in\bC^N$ and $p'_0\in\bC^{N'}$,
respectively, such that $M$ is of finite
type at $p_0$ and $M'$ is
$\ell_0$-nondegenerate at $p'_0$ for some
integer
$\ell_0\geq 0$. Let $d$ be the codimension of $M$.
Then
there exist a finite number of formal power
series mappings of the form
$$
\Psi^k(Z,\Lambda)\sim\sum_{\alpha}
\frac{Q^k_\alpha(\Lambda)}
{P^k(\Lambda)^{l^k_\alpha}}(Z-p_0)^\alpha,\quad
k=1\ldots, l, \tag0.2
$$
where $P^k$ and $Q^k_\alpha$ are ($\bC$ and
$\bC^{N'}$ valued, respectively) polynomials
on the jet space
$J^{(d+1)\ell_0}(\bC^N,\bC^{N'})_
{(p_0,p'_0)}$ and $l^k_\alpha$ are nonnegative integers, such 
that the following 
holds. For any smooth CR submersive mapping
$H\:(M,p_0)\to (M',p'_0)$, and more generally any formal 
(holomorphic) 
CR submersive power series mapping $\hat H\:(M,p_0)\to 
(M',p_0')$, there exists $k\in\{1,\ldots,l\}$ with $P^{k} 
\left(j_{p_0}^{(d+1)\ell_0}(\hat
H)\right)\neq 0$ if $d+1$ is even and
$P^{k} 
\left(\overline{j_{p_0}^{(d+1)\ell_0}(\hat
H)}\right)\neq 0$ if $d+1$ is odd, and for any such $k$,
$$
\aligned
\hat H(Z) &\sim
\Psi^k\left (Z,j_{p_0}^{(d+1)\ell_0}(\hat
H)\right),\ \text{\rm  if $d+1$ is even},\\
\hat H(Z) &\sim
\Psi^k\left(Z,\overline{j_{p_0}
^{(d+1)\ell_0}(\hat H)}\right),\  
\text{\rm  if $d+1$ is odd}.
\endaligned\tag0.3
$$
If, in addition, $M$ and $M'$ are
real-analytic, then the series \thetag{0.2},
for $k\in\{1,\ldots, l\}$, converges in a
neighborhood of
$(p_0,\Lambda_0)$ in $\bC^N\times
J^{(d+1)\ell_0}(\bC^N,\bC^{N'})_
{(p_0,p'_0)}$ for every $\Lambda_0$
satisfying $P^k(\Lambda_0)\neq 0$.
\endproclaim

Theorem 1, which will follow from the more
general Theorem 2.1.5, has a number of
applications.  Our first application, which
will be given in a more general
form in Theorem 2.1.1, states that a
holomorphic mapping sending
$M$ into
$M'$ which is CR submersive at $p_0$ is
uniquely determined by finitely many
derivatives at
$p_0$.

\proclaim{Theorem 2}  Let $M$, $p_0$, $d$,
$M'$,
$p_0'$, and $\ell_0$ be as in Theorem 1.
Then there exists an integer $k_0$,
depending only on $M$, with
$1<k_0\leq d+1$ such that the following
holds. If 
$H^1,H^2\:(\bC^N,p_0)\to (\bC^{N'}, p'_0)$
are holomorphic mappings near $p_0$ such
that $H^j(M)\subset M'$, $H^j$ is CR
submersive at $p_0$, for $j=1,2$,
and 
$$
\frac{\partial^{|\alpha|}H^1}{\partial
Z^\alpha}(p_0)=
\frac{\partial^{|\alpha|}H^2}{\partial
Z^\alpha}(p_0),\quad\forall
\alpha\:|\alpha|\leq k_0\ell_0,\tag0.4
$$
then $H^1\equiv H^2$.
\endproclaim

The conditions of finite type and finite nondegeneracy in Theorem 2 are
also essentially necessary, in a certain sense, for the conclusion to
hold. We refer the reader to the discussion in \S2.2. 

Our second application, Theorem 3 below
(which is an easy consequence of Theorem 1),
deals with real-analytic submanifolds. It
gives sufficient conditions for all CR
submersive formal mappings between
real-analytic generic submanifolds to be
convergent. 

\proclaim{Theorem 3} Let $M$ and $M'$
be real-analytic generic submanifolds
through
$p_0\in\bC^N$ and $p'_0\in\bC^{N'}$,
respectively, such that $M$ is of finite
type at $p_0$ and $M'$ is
finitely nondegenerate at $p'_0$. Then, any formal  
(holomorphic) CR submersive
mapping $\hat H\:(M,p_0)\to
(M,p_0')$ is
convergent, i.e. $\hat H$ is the Taylor
series at $p_0$ of a holomorphic mapping
$H\:(\bC^N,p_0)\to (\bC^{N'},p_0')$ near
$p_0$ with
$H(M)\subset M'$.
In particular, 
$(M,p_0)$ and $(M',p_0')$ are formally
equivalent if and only if they are
biholomorphically equivalent.
\endproclaim

For our next application of Theorem 1, we shall denote
the set of holomorphic mappings
$H\:(\bC^N,p_0)\to (\bC^{N'},p_0')$ which
map $M$ into $M'$ and are CR submersive at
$p_0$ by $\Cal F(M,p_0;M',p_0')$.
This set has a natural inductive limit
topology induced by uniform convergence on
compact neighborhoods of
$p_0$. We have the following result, which
will be proved in \S4.3.

\proclaim{Theorem 4} Let $M$, $p_0$, $M'$,
and $p_0'$ be as in Theorem $3$.
Let $d$ be the codimension of $M$, and $\ell_0$ a nonnegative integer
such that $M'$ is $\ell_0$-nondegenerate at $p_0'$. Then 
there exist two real
algebraic subvarieties
$$
A, B\subset
J^{(d+1)\ell_0}(\bC^{N},\bC^{N'})_
{(p_0,p'_0)}\tag 0.5
$$ 
such that the 
mapping
$$
j^{(d+1)\ell_0}_0\:\Cal F(M,p_0;M',p_0')\to
J^{(d+1)\ell_0}(\bC^{N},\bC^{N'})_{(p_0,p'_0)}
\tag0.6
$$
is a homeomorphism onto 
$A\setminus B$. In addition,
the image $A\setminus B$ is totally real at
each nonsingular point.   
\endproclaim

The case where $(M,p_0)=(M',p'_0)$ is of
particular interest. In this case, the
set of mappings $\Cal F(M,p_0;M,p_0)$
consists of biholomorphisms of
$(\bC^N,p_0)$ (see Corollary 1.27)
and, hence, forms a group under composition.
This group is called the {\it stability
group} of $M$ at $p_0$, and is denoted
$\Aut(M,p_0)$. It follows from Theorem 4
that $\Aut (M,p_0)$, where $M$ is a real-analytic generic submanifold
of codimension $d$
which is $\ell_0$-nondegenerate and of finite type at $p_0\in M$, is a
real Lie group
which can be homeomorphically embedded as
an algebraic, totally real subgroup of the
jet group $G^{(d+1)\ell_0}(\bC^N)_{p_0}$ (see also
Theorem 2.1.14).  Here $G^{(d+1)\ell_0}(\bC^N)_{p_0}$
consists of those elements in
$J^{(d+1)\ell_0}(\bC^{N},\bC^{N})_ {(p_0,p_0)}$ which
are invertible.

Recall that a real-analytic generic
submanifold $M\subset\bC^N$ is called  {\it
holomorphically nondegenerate} at $p\in 
M$ if there are no non-trivial holomorphic vector 
fields (i.e. holomorphic sections of the holomorphic
tangent bundle 
$T'\bC^N$) near $p$
which are tangent to $M$.  If $M$ is connected, then
it is either holomorphically nondegenerate at every
point or at no point (see e.g.\ [BER4, Theorem
11.5.1]). We say that a connected real-analytic
generic submanifold $M$ is holomorphically
nondegenerate if it is so at some (and hence at every)
point. The relation between  holomorphic nondegeneracy
and finite nondegeneracy is discussed in \S2.2. The
following result is then a corollary of  Theorem 4 and
the discussion in \S2.2. 

\proclaim{Theorem 5} Let $M$ be a connected,
real-analytic, holomorphically nondegenerate, generic submanifold of
codimension $d$ in 
$\bC^N$ which is of finite type at some point. Then there
exists a proper real-analytic subvariety $V$ of $M$ such that the
following holds for every $p\in M\setminus V$. The jet mapping 
$$
j^{(d+1)(N-d)}_p\:\Aut(M,p)\to
G^{(d+1)(N-d)}(\bC^N)_p
$$
is a continuous injective group
homomorphism which is a homeomorphism onto a totally
real algebraic Lie subgroup of
$G^{(d+1)(N-d)}(\bC^N)_p$.
\endproclaim

In \S5, we consider smooth perturbations of generic submanifolds,
satisfying the appropriate conditions, and
study the behavior of the functions $\Psi^j$ in Theorem
1 under such perturbations (Theorem 5.1.1). As a
consequence (Theorem 5.1.9) we obtain the result that if
the stability group of a real-analytic generic
submanifold $M$ is discrete, then it remains discrete
under real-analytic (small) perturbations of
$M$. One of the more important examples of a perturbation is allowing the
base point $p_0\in M$ to vary. As another application of Theorem
5.1.1, we show that the
topological space $\cup_{p\in M,p'\in M'} \Cal F(M,p;M',p')$ is
homeomorphically embedded in the jet manifold
$J^{(d+1)\ell_0}(\bC^N,\bC^{N'})$ (see \S5.2) as a difference
$A\setminus B$, where $A,B\subset J^{(d+1)\ell_0}(\bC^N,\bC^{N'})$ 
are real-analytic subvarieties whose fibers in
$J^{(d+1)\ell_0}(\bC^N,\bC^{N'})_{(p,p')}$ are real-algebraic (see
Theorem 5.2.9). 

The paper concludes with an application, in \S6, of our
methods to the study of algebraicity of holomorphic
mappings which map one real algebraic submanifold into
another; the reader is referred to  Theorem 6.1 for this
result. 

An important tool in the proofs of the
results in this paper is the sequence of Segre mappings
(see \S3.1), which in the
real-analytic case were introduced in
[BER1] along with the Segre sets. The Segre variety of a real-analytic
hypersurface
$M$, which coincides with the first Segre
set, was first introduced by Segre [Se]. Its
use in the study of holomorphic mappings
between real-analytic hypersurfaces was
pioneered by Webster [W1], [W3]. Since then,
its use has been crucial in the work of 
many mathematicians, including
Diederich--Webster [DW],
Diederich--Fornaess [DF], Forstneric [F],
Huang [Hu], Diederich--Pinchuk [DP], and
others. (See also the notes in [BER4,
Chapter X--XII].) The Segre sets,
introduced in [BER1] and playing a crucial
role in the proofs in  [BER1] and [BER2],
have also been an important tool in
the work of Zaitsev [Z2], [Z3]. We should also
mention here the work of Christ, Nagel,
Stein, and Wainger [CNSW] in a different
context, in which they study the relation
between certain curvature conditions on
families of submanifolds in $\bR^n$. Two of
these conditions are strikingly similar to
the two equivalent conditions in Theorem
3.1.9. An alternative presentation of the Segre sets and mappings,
along with a more detailed study, will be given in
the forthcoming paper [BER5].

The study of automorphism groups of
bounded domains in $\bC^N$ goes back to 
H. Cartan [HC] (and was later continued by
Kaneyuki [Ka] and, more recently, by Zaitsev
[Z1]). The structure of the local
transformation groups of Levi nondegenerate
hypersurfaces in $\bC^2$ was investigated
by E. Cartan [EC1], [EC2] in connection with his
work on biholomorphic equivalence problem. 
His results were later
extended to Levi nondegenerate
hypersurfaces in higher dimensions by
Tanaka [Ta] and Chern--Moser [CM]. 
In particular, the conclusions of Theorems 2,3,
and 4 for real-analytic Levi
nondegenerate hypersurfaces follow from
their work. The convergence of the formal series 
\thetag{0.2} in Theorem 1 in the real-analytic case
seems to be new even for Levi nondegenerate
hypersurfaces. 
Further results on
transformation groups of Levi
nondegenerate hypersurfaces were obtained
by a number of mathematicians,
including Webster [W2], Burns-Shnider [BS],
and the Moscow school (Beloshapka,
Krushilin, Loboda, Vitushkin, etc.; see
Krushilin [Kr], and
Vitushkin [Vi]). Stanton [St1], [St2]
considered infinitesimal CR automorphisms
on general real-analytic hypersurfaces.
(See also [BER2] for results on
infinitesimal CR automorphism in higher
codimensions.) The case of higher
codimensions was considered by
Tumanov--Henkin [TH], Tumanov [Tu] in the
case of quadratic manifolds, and by Beloshapka [B] in the more general
case where the Levi forms of the submanifolds are nondegenerate. For
these classes
of manifolds, the conclusion of Theorem 2 follows from their
work. 

We conclude the introduction by giving a
brief history of results related to
those of Theorems 1--4 above.  Theorem 2, in
the case
$N=N'$ with
$M$ and $M'$ real-analytic and of the same
dimension, was obtained in [BER2]. In
[BER3], Theorems 1, 3, and 4, with slightly
weaker conclusions, were proved in the case
$N=N'$ with
$M$ and
$M'$ real-analytic hypersurfaces. Zaitsev
[Z2] proved a weaker version of Theorem 1
for real-analytic 
CR submersive mappings
(in the real-analytic case), namely one in
which the jet space of order
$(d+1)\ell_0$ is replaced by that of order 
$2(d+1)\ell_0$ and where the dependence on
the jets is only local analytic instead of
rational. In particular, his result shows
that the stability group is a Lie group
with the natural topology. However, for the application 
given by
Theorem 3 it is 
important to prove Theorem 1 for formal CR submersive mappings.

For most of the proofs of the results
mentioned above, it is convenient to work
with formal mappings between formal generic
submanifolds. Hence, most results presented
here will be reformulated, and proved, in this
more general context. The following section presents the necessary
preliminaries and definitions. In what follows, the distinguished
points $p_0$ and $p_0'$ on $M$ and $M'$, respectively, will, for convenience 
and without loss of
generality, be assumed to be $0$.

\heading \S1. Preliminaries on formal submanifolds
and mappings\endheading

Let
$\bC[[x]]=\bC[[x_1,\ldots,x_k]]$ be the ring of
formal power series in $x=(x_1,\ldots, x_k)$
with complex coefficients. Suppose that
$\rho=(\rho_1,\ldots,
\rho_d)\in\bC[[Z,\zeta]]^d$, where
$Z=(Z_1,\ldots, Z_N)$ and
$\zeta=(\zeta_1,\ldots,\zeta_N)$, satisfies the
reality condition 
$$
\rho(Z,\zeta)\sim\bar\rho(\zeta,Z),\tag1.1
$$
where $\bar\rho$ is the formal series
obtained from $\rho$ by replacing each
coefficient in the series by its complex
conjugate; we use the symbol $\sim$ to denote equality of formal
power series.
If, in addition, the series
$\rho$ satisfies the condition
$\rho(0)=0$, and
$$
\partial\rho_1(0)\wedge
\ldots\wedge \partial\rho_d(0)\neq 0,\tag1.2
$$
then we say that $\rho$ defines a {\it formal
real submanifold} $M$ of $\bC^N$ through $0$
of codimension $d$ (and dimension $2N-d$). 
If $M'$
is another such formal real submanifold defined
by
$\rho'=(\rho'_1,\ldots,\rho'_d)$, then we
shall say that
$M=M'$ if there exists a $d\times d$ matrix of
formal power series
$a(Z,\zeta)$ (necessarily invertible at
$0$) such that
$$
\rho(Z,\zeta)\sim
a(Z,\zeta)\rho'(Z,\zeta).\tag1.3
$$ 

These definitions are
motivated by the fact that if in addition the
components of
$\rho$ are convergent power series, then the
equations
$\rho(Z,\bar Z)=0$ define a real-analytic
submanifold $M$ of $\bC^N$ through $0$.
Moreover, if $M'$ is another such defined by a
convergent power series $\rho'$, then
\thetag{1.3} holds if and only if $M$ and
$M'$ are the same. Also, if $M$ is a smooth real
submanifold in $\bC^N$ through $0$, then the
Taylor series at $0$ of a smooth defining
function $\rho(Z,\bar Z)$ of $M$ near $0$, with $\bar Z$ formally
replaced by 
$\zeta$,  defines a formal  
real submanifold through $0$. These observations will be used to
deduce the results given in the introduction from the corresponding results
for formal real submanifolds.

If the formal series $\rho$ defining $M$
satisfies the stronger condition
$$
\partial_Z\rho_1(0)\wedge
\ldots\wedge\partial_Z\rho_d(0)\neq 0,\tag1.4
$$
(which in particular implies \thetag{1.2})
then we say that the formal real submanifold
$M$ is {\it generic}. 
\comment
More generally, we say
that $M$ is {\it CR} if
$$
\rk\,(\partial_Z\rho_1(0),\ldots,
\partial_Z\rho_d(0))=
\rk\,(\partial_Z\rho_1(Z,\zeta),\ldots,
\partial_Z\rho_d(Z,\zeta));\tag1.5 
$$
here, the right hand side of \thetag{1.5} denotes the
rank of the matrix
$$(\partial_Z\rho_1(Z,\zeta),\ldots,
\partial_Z\rho_d(Z,\zeta)),$$ where the
 $\partial_Z\rho_j(Z,\zeta)$, $j=1,\ldots, d$, are 
considered as vectors in $\Bbb K^N$, with $\Bbb K$ 
the quotient field of $\bC[[Z,\zeta]]$.    
\endcomment
We say that a formal vector field 
$$
X=\sum_{j=1}^N\left(a_j(Z,\zeta)\frac{\partial}
{\partial Z_j}+b_j(Z,\zeta)\frac{\partial}
{\partial \zeta_j}\right),\tag1.5
$$
with $a_j,b_j\in\bC[[Z,\zeta]]$,
is tangent to the formal real submanifold $M$
if 
$$
X\rho(Z,\zeta)\sim
a(Z,\zeta)\rho(Z,\zeta),\tag1.6
$$ 
for some $d\times d$ matrix of formal power series $a(Z,\zeta)$.

We say
that the formal vector field $X$ in
\thetag{1.6} is of type
$(0,1)$ if $a_j\sim0$, $j=1,\ldots, N$, and
similarly of type $(1,0)$ if the $b_j\sim
0$. Let $\Cal D_M$ denote the
$\bC[[Z,\zeta]]$-module generated by all formal $(0,1)$ and
$(1,0)$ vector fields tangent to $M$, and $\frak
g_M$ the Lie algebra generated by $\Cal D_M$. We
denote by
$\frak g_M(0)$ the complex 
vector space obtained by evaluating the
coefficients of the formal vector fields in
$\frak g_M$ at $0$. Similarly, we use the notation $\Cal D_M(0)$ for
the complex  
vector space obtained by evaluating the
coefficients of the formal vector fields in
$\Cal D_M$ at $0$. (The reader should observe the analogy, for a smooth
real submanifold $M$ through $0$, between the
complexified complex tangent space $\bC T^c_0 M$ to $M$ at $0$ and
$\Cal D_M(0)$ for the corresponding formal submanifold.) Thus, we have
$\Cal D_M(0)\subset \frak g_M(0)\subset T'_0\bC^{2N}$, where
$T'_0\bC^{2N}$ denotes the holomorphic
tangent space of $\bC^{2N}$ at $0$,

We say that $M$ is of {\it
finite type} at
$0$ if $\dim_\bC \frak g_M(0)=\dim M=2N-d$. 
(Note that the vector space $\Cal D_M(0)$ has
dimension $2N-2d$;  this follows easily from 
the fact that $M$ is generic and of
codimension
$d$.) 

We shall also need the notion of finite
nondegeneracy of a formal generic submanifold.
We say that the formal generic submanifold $M$
is {\it finitely nondegenerate} at $0$ if there
exists an integer $\ell\geq 0$ such that 
$$
\Span
\left\{L^\alpha\left(\frac{\partial\rho_j}
{\partial Z}\right)(0)\:1\leq j\leq
d,\ |\alpha|\leq \ell\right\}=\bC^N.\tag1.7
$$
Here, $L_1,\ldots, L_n$ is a basis for the
$\bC[[Z,\zeta]]$-module of all formal $(0,1)$ vector
fields tangent to $M$ (so
$n=N-d$) modulo those whose coefficients are in the
ideal generated by $\rho_1,\ldots,\rho_d$. We  also use
multi-index notation, i.e.\ we introduce the vector 
$L=(L_1,\ldots, L_n)$
and, for any
$\alpha\in\Bbb Z_+^n$, we
write 
$$
L^\alpha=L_1^{\alpha_1}\ldots
L_n^{\alpha_n},\quad |\alpha|=
\sum_{j=1}^n\alpha_j.\tag1.8
$$ 
More precisely, we say
that $M$ is {\it $\ell_0$-nondegenerate} at $0$ if
$\ell_0$ is the smallest integer for which
\thetag{1.7} holds. It is an easy exercise to
show that the definition of
$\ell_0$-nondegeneracy (and hence that of finite
nondegeneracy) does not depend on the choice of
basis
$L=(L_1,\ldots, L_n)$,  defining series
$\rho=(\rho_1,\ldots,\rho_d)$, or the choice of
coordinates $Z$. Hence, $\ell_0$-nondegeneracy is
a property of the formal generic submanifold $M$. The
reader is also referred to [BER4] for further
discussion of these notions, as well as that of
finite type, for smooth and
real-analytic generic submanifolds.

Let
$H\:(\bC^N,0)\to (\bC^{N'},0)$ be a formal
mapping, i.e. $H\in \bC[[Z_1,\ldots, Z_N]]^{N'}$ such that each
component of $H(Z)=(H_1(Z),\ldots, H_{N'}(Z))$ has no constant term. To such a
formal mapping $H$ we 
associate a formal mapping
$\Cal H\:(\bC^N\times\bC^N,0)\to
(\bC^{N'}\times\bC^{N'},0)$ defined by
$$
\Cal H(Z,\zeta)\sim (H(Z),\bar
H(\zeta)).\tag1.9
$$ 
If $M$ and $M'$ are formal real
submanifolds of
$\bC^N$ and $\bC^{N'}$ defined by formal
series
$\rho(Z,\zeta)=(\rho_1(Z,\zeta),\ldots,
\rho_d(Z,\zeta))$ and
$\rho'(Z',\zeta')=(\rho'_1(Z',\zeta'),\ldots,
\rho'_{d'}(Z',\zeta'))$, respectively,
then we say
that the formal mapping $H$, as above, maps
$M$ into $M'$, denoted $H\:(M,0)\to
(M',0)$, if
$$
\rho'(H(Z),\bar H(\zeta))\sim
c(Z,\zeta)\rho(Z,\zeta),\tag1.10
$$ 
for some $d'\times d$ 
matrix $c(Z,\zeta)$ of formal power series. It will be
be convenient to choose {\it normal
coordinates}, $Z=(z,w)$ and
$\zeta=(\chi,\tau)$
with
$z=(z_1,\ldots, z_n)$, $w=(w_1,\ldots,
w_d)$ (so
$n+d=N)$, $\chi=(\chi_1,\ldots,\chi_n)$,
and $\tau=(\tau_1,\ldots,\tau_d)$, in $\bC^N\times\bC^N$ for $M$ at $0$. By this we mean
 a formal change of coordinates $Z=Z(z,w)$ with $Z(z,w)$
 a formal invertible mapping $(\bC^N,0)\to (\bC^N,0)$,
and $\zeta=\bar Z(\chi,\tau)$ the corresponding change,
 such that 
$$
\rho(Z(z,w),\bar Z(\chi,\tau))\sim a(z,w,\chi,\tau)
(w-Q(z,\chi,\tau)),\tag1.11
$$
where $a(z,w,\chi,\tau)$ is an invertible $d\times d$ 
matrix of formal power series, and the vector valued 
$Q\in\bC[[z,\chi,\tau]]^d$ satisfies 
$$
Q_j(0,\chi,\tau)\sim
Q_j(z,0,\tau)\sim\tau_j, \quad j=1,\ldots,d.\tag1.12
$$
(See [BER4, Chapter IV.2].) It follows
from the reality of $\rho$ that
$$
\rho(Z(z,w),\bar Z(\chi,\tau))\sim b(z,w,\chi,\tau)(\tau-\bar Q(\chi,z,w)),\tag1.13
$$
where 
$b(z,w,\chi,\tau)$ is an invertible $d\times d$ matrix of formal
power series. Similarly, we choose 
normal coordinates
$Z'=(z',w')$ and
$\zeta'=(\chi',\tau')$ in
$\bC^{N'}\times\bC^{N'}$, with
$z'=(z'_1,\ldots, z'_{n'})$ and
$w'=(w'_1,\ldots, w'_{d'})$ (so
$n'+d'=N')$, such that $M'$ is
defined by $w'-Q'(z',\chi',\tau')$ (or, more precisely,
by $a'(z,w,\chi,\tau)(w'-Q'(z',\chi',\tau'))$, for some matrix $a'$
making the expression real). 
Then we may write the formal mapping
$H=(F,G)$, with $F=(F_1,\ldots, F_{n'})$
and $G=(G_1,\ldots, G_{d'})$, and the
condition $H\:(M,0)\to
(M',0)$ can be expressed by either of the
equations
$$
\aligned
G(z,w) \sim Q'(F(z,w&),\bar
F(\chi,\tau),\bar G(\chi,\tau))\\
\text{\rm or}&\\ \bar G(\chi,\tau) \sim
\bar Q'(\bar F(\chi,\tau)&,
F(z,w),
G(z,w)),\endaligned\tag1.14
$$
for $\tau= \bar Q(\chi,z,w)$ or
$w= Q(z,\chi,\tau)$. An observation
that will be useful is that
$G(z,0)\sim 0$, as is easily verified by
taking $w=\tau=0$ and $\chi=0$ in
\thetag{1.14}.

Note that if $M$ and $M'$ correspond to
real-analytic submanifolds and the formal
mapping
$H$ defines a holomorphic mapping
$(\bC^N,0)\to (\bC^{N'},0)$ in some
neighborhood of the origin, then
$H\:(M,0)\to (M',0)$ if and only if
$H(M)\subset M'$. Moreover, if
$M\subset\bC^N$ and
$M'\subset\bC^{N'}$ are smooth CR
submanifolds through the origin and
$h\:M\to M'$ is a smooth CR mapping, with
$h(0)=0$, then there exists a unique
formal mapping
$H\:(M,0)\to (M',0)$ such that, for any
local parametrization 
$$
\bR^{2N-d}\supset U\ni x\mapsto Z(x)\in M,\tag
1.15
$$ with $Z(0)=0$, we have 
$$ h(Z(x))\sim H(Z(x)),\tag1.16
$$ where the left hand side of
\thetag{1.16} refers to the Taylor
expansion at
$0$ of the smooth mapping $x\mapsto
h(Z(x))$, and the right hand side is taken
in the sense of composition of $H$ and the
Taylor series at $0$ of $Z(x)$. (See e.g.\
[BER4,
\S1.7]).  Observe that if $h$ is the restriction to $M$ of a
holomorphic mapping $(\bC^N,0)\to (\bC^{N'},0)$, then $H(Z)$ is the
Taylor series of this holomorphic mapping.

Given a formal mapping $H\:(\bC^N,0)\to
(\bC^{N'},0)$, we denote the
tangent mapping of the associated formal
mapping
$\Cal H\: (\bC^N\times\bC^N,0)\to
(\bC^{N'}\times\bC^{N'},0)$  by $d\Cal
H\:T'_0\bC^{2N}\to T'_0\bC^{2N'}$; thus,
$d\Cal H$ is the mapping taking a vector
$$
X=\sum_{j=1}^N\left(a_j\frac{\partial}{\partial
Z_j}+ b_j\frac{\partial}{\partial
\zeta_j}\right),\quad
a_j,b_j\in\bC\tag1.17
$$  in $T'_0\bC^N$ to the vector
$$  
d\Cal
H(X)=\sum_{j=1}^{N'}\left((XH_j(Z))(0)\frac
{\partial}{\partial Z'_j}+ (X\bar
H_j(\zeta))(0)\frac{\partial}{\partial
\zeta'_j}\right)\tag1.18
$$ 
in $T'_0\bC^{N'}$. Recall
that we say that the vector $X$ in
\thetag{1.17} is a
$(1,0)$ vector if $b_j=0$, $j=1,\ldots,
N$, and
$(0,1)$ vector if $a_j=0$, $j=1,\ldots,
N$. It is clear that $d\Cal H$ maps
$(0,1)$ vectors to
$(0,1)$ vectors and $(1,0)$
vectors to $(1,0)$ vectors. If $H\:(M,0)\to
(M',0)$, then it follows that $d\Cal
H$ maps $\Cal D_M(0)$ into $\Cal
D_{M'}(0)$. We shall say that a formal mapping $H\:(M,0)\to (M',0)$ is
{\it CR submersive} at $0$ if 
$$
d\Cal H(\Cal
D_M(0))=\Cal D_{M'}(0). \tag1.19
$$

\proclaim{Proposition 1.20} Let $M$ and
$M'$ be formal generic submanifolds
through the origin in $\bC^N$ and
$\bC^{N'}$, respectively. If
$H\:(M,0)\to (M',0)$ is a CR submersive formal mapping,
then $d\Cal
H(\frak g_M(0))=\frak g_{M'}(0)$.
\endproclaim

\demo{Proof} Let $L'_1,\ldots, L'_{n'}$ be
a basis for the formal $(1,0)$ vector fields
tangent to $M'$, and $\tilde L'_1,\ldots,
\tilde L'_{n'}$ a basis for the formal $(0,1)$
vector fields tangent to $M'$. We claim
that for any formal $(1,0)$ vector field $L$
tangent to $M$ there exist formal power
series $a_j=a_j(Z,\zeta)$, $j=1,\ldots n'$,
such that for any $f\in\bC[[Z',\zeta']]$
we have
$$
L(f\circ \Cal H)\sim\sum_{j=1}^{n'}
a_j((L'_jf)\circ\Cal H)\tag1.21
$$
as power series in $(Z,\zeta)$. To prove
the claim, it suffices to find 
$a_j(Z,\zeta)$ satisfying
\thetag{1.21} with $f(Z',\zeta')= Z'_k$,
$k=1,\ldots N'$. This is done by using the
chain rule and elementary linear algebra
over the quotient field
of $\bC[[Z,\zeta]]$. The details are
left to the reader. Similarly, for any
$(0,1)$ vector field 
$\tilde L$ tangent to $M$, we can find
$\tilde a_j(Z,\zeta)$ such that
$$
\tilde L(f\circ \Cal H)\sim\sum_{j=1}^{n'}
\tilde a_j((\tilde L'_jf)\circ\Cal
H)\tag1.22
$$
as power series in $(Z,\zeta)$. Thus, for
bases $L_1,\ldots, L_n$ of the $(1,0)$
vector fields tangent to $M$ and $\tilde
L_1,\ldots, \tilde L_n$ of the $(0,1)$
vector fields tangent to $M$, we obtain
two $n\times n'$ matrices
$(a_{jk}(Z,\zeta))$ and $(\tilde
a_{jk}(Z,\zeta))$ of formal power series
such that
$$
L_j(f\circ \Cal H)\sim\sum_{k=1}^{n'}
a_{jk}((L'_kf)\circ\Cal H),\quad 
\tilde L_j(f\circ \Cal
H)\sim\sum_{k=1}^{n'}
\tilde a_{jk}((\tilde L'_kf)\circ\Cal
H),\tag1.23
$$
for all $f\in\bC[[Z',\zeta']]$. It is easy
to verify that \thetag{1.19} implies that
the rank of each of these matrices at
$0$ equals $n'$. Hence, we may assume,
after a linear transformation of the
$L_j$'s and
$\tilde L_j$'s (over
$\bC[[Z,\zeta]]$) if necessary,  that  
$$
L_j(f\circ \Cal H)\sim
(L'_jf)\circ\Cal H,\quad 
\tilde L_j(f\circ \Cal
H)\sim (\tilde L'_jf)\circ\Cal
H,\quad j=1,\ldots, n',\tag1.24
$$
for all $f\in\bC[[Z',\zeta']]$. It follows
immediately from \thetag{1.24} that we
also have
$$
[X,Y](f\circ \Cal H)\sim
([X',Y']f)\circ\Cal H,\tag1.25
$$
for any $X,Y\in\{L_1,\ldots, L_{n'},\tilde
L_1,\ldots, \tilde L_{n'}\}$ and
corresponding
$X',Y'\in\{L'_1,\ldots,$ $L'_{n'},
\tilde L'_1,\ldots, \tilde L'_{n'}\}$
(i.e. such that $X$, $X'$ and $Y$,
$Y'$ satisfy 
\thetag{1.24}). In particular,  we have 
$$
d\Cal H([X,Y]_0)=[X',Y']_0.\tag1.26
$$
Repeating this argument for commutators of
any length, we can conclude that
$d\Cal H(\frak g_M(0))=\frak g_{M'}(0)$.
This completes the proof of Proposition
1.20.
\qed\enddemo

\proclaim{Corollary 1.27} Let $M$ and 
$M'$ be formal generic submanifolds of
codimension $d$ and $d'$ through the origin
in $\bC^N$ and $\bC^{N'}$, respectively.
Assume that $M'$ is of finite type at $0$,
and let $H\:(M,0)\to (M',0)$ be a CR submersive formal
mapping. Then $d
H(T'_0\bC^{N})=T'_0\bC^{N'}$ and $d\geq
d'$. If in addition $N=N'$, then $d
H$ is an isomorphism of $T'_0\bC^{N}$ into
itself, i.e. the formal mapping
$H$ is invertible. 
\endproclaim
\demo{Proof} Since $d\Cal H$ maps $\Cal
D_M(0)$ onto $\Cal D_{M'}(0)$ by
assumption, and hence $\frak g_M(0)$ onto
$\frak g_{M'}(0)$ by Proposition 1.20, it
follows that the induced mapping from
the quotient space $\frak g_M(0)/\Cal
D_M(0)$ is onto $\frak g_{M'}(0)/\Cal
D_{M'}(0)$. Since $M'$ is of finite type
at $0$, it follows that $\dim \frak
g_{M'}(0)/\Cal D_{M'}(0)=d'$. Since $d\geq
\dim \frak
g_{M}(0)/\Cal D_{M}(0)$, it follows that
$d\geq d'$. (Hence, we also have $N\geq
N'$.)

To prove that
$d H(T'_0\bC^{N})=T'_0\bC^{N'}$, we take
normal coordinates $Z=(z,w)$,
$\zeta=(\chi,\tau)$ in $\bC^{2N}$ for $M$
and  $Z'=(z',w')$,
$\zeta'=(\chi',\tau')$ in $\bC^{2N'}$ for
$M'$. For $H=(F,G)$, the fact that
$H\:(M,0)\to (M',0)$ is expressed by
\thetag{1.14}. Also, $d \Cal H(\Cal
D_M(0))=\Cal D_{M'}(0)$ is equivalent to
$\partial F/\partial z(0)$ having rank
$n'$. Using the facts that
$G(z,0)\sim 0$ and $M'$ is of finite type
at $0$, and applying Proposition 1.20, we
conclude that the rank of $\partial
G/\partial w(0)$ is $d'$. This completes
the proof of $d
H(T'_0\bC^{N})=T'_0\bC^{N'}$. The second
statement of Corollary 1.27 is an
immediate consequence of the first. \qed 
\enddemo

\heading \S 2. Uniqueness and
parametrization of formal
mappings\endheading

\subhead \S2.1. Main results\endsubhead In this section, we shall give
results on uniqueness and 
parametrization of formal 
mappings between formal real submanifolds from which Theorems 1--3
presented in the introduction will follow.  We first give 
sufficient conditions so that a mapping
sending $M$ into $M'$ is determined by a
finite number of derivatives of the mapping
at $0$. The necessity of these conditions will be discussed in \S2.2.

\proclaim{Theorem 2.1.1} Let $M$ and $M'$
be formal generic submanifolds through
$0\in\bC^N$ and $0\in\bC^{N'}$,
respectively, such that $M$ is of finite
type at $0$ and $M'$ is
$\ell_0$-nondegenerate at $0$ for some
integer
$\ell_0$. Let $d$ denote the codimension of $M$. Then there exists an
integer $k_0$, depending only on $M$, with $1<k_0\leq d+1$, 
such that the following holds. If
$H^1,H^2\:(M,0)\to (M',0)$ are CR submersive formal
mappings such that 
$$
\frac{\partial^{|\alpha|}H^1}{\partial
Z^\alpha}(0)=
\frac{\partial^{|\alpha|}H^2}{\partial
Z^\alpha}(0),\quad\forall
\alpha\:|\alpha|\leq k_0\ell_0,\tag2.1.2
$$ 
then
$H^1\sim H^2$.
\endproclaim

\remark{Remark $2.1.3$} If $N=N'$, $\dim M=\dim M'$, and
$H^j$,
$j=1,2$, are invertible formal mappings,
then 
$d\Cal H^j:\Cal D_M(0)\to\Cal D_{M'}(0)$,
$j=1,2$, are necessarily isomorphisms, and
hence surjective. More generally, if 
$n=n'$ (recall that $n$ denotes $N-d$,
where
$d$ denotes the codimension of $M$, and
similarly for $n'$ and $M'$) and the
mappings
$d\Cal H^j:\Cal D_M(0)\to\Cal D_{M'}(0)$,
$j=1,2$, are injective,
then they are also necessarily
surjective. (Indeed, $\dim \Cal
D_M(0)=2n$.) \endremark\medskip

It is clear from the remarks in \S1 that Theorem 2.1.1 is a more general
version of Theorem 2 in the introduction. The proof of Theorem 2.1.1
will be given in \S3.4.

Let $E(\bC^N,\bC^{N'})_{(0,0)}$ denote the
set of germs of holomorphic mappings
$(\bC^N,0)\to (\bC^{N'},0)$ and $\hat
E(\bC^N,\bC^{N'})_{(0,0)}$ the set of
formal mappings $(\bC^N,0)\to
(\bC^{N'},0)$. For each positive integer
$k$, we denote by
$J^k(\bC^N,\bC^{N'})_{(0,0)}$ the jet
space of order $k$ of holomorphic mappings 
$(\bC^N,0)\to (\bC^{N'},0)$, and by
$j^k_0\:\hat E(\bC^N,\bC^{N'})_{(0,0)}\to
J^k(\bC^N,\bC^{N'})_{(0,0)}$ the jet
mapping taking a formal mapping $H$ to its
$k$th jet at $0$,
$j^k_0 (H)$. In particular,
$J^1(\bC^N,\bC^{N'})_{(0,0)}$ can be
viewed as the space of linear mappings from
$\bC^N$ to $\bC^{N'}$. For $k\geq l\geq 
1$, we denote by
$$
j_0^{k,l}\:J^{k}(\bC^N,\bC^{N'})_{(0,0)}\to
J^l(\bC^N,\bC^{N'})_{(0,0)}\tag2.1.4
$$ 
the canonical mapping induced by 
$j^l_0\:\hat E(\bC^N,\bC^{N'})_{(0,0)}\to
J^l(\bC^N,\bC^{N'})_{(0,0)}$.

Given coordinates
$Z$ and
$Z'$ on
$\bC^N$ and $\bC^{N'}$, the jet space
$J^k(\bC^N,\bC^{N'})_{(0,0)}$ can be
identified with the set of polynomial
mappings $(\bC^N,0)\to (\bC^{N'},0)$ of
degree $k$. The coordinates on $J^k(\bC^N,\bC^{N'})_{(0,0)}$, 
which we will
denote by $\Lambda$, can then be taken to
be the coefficients of these polynomials.
Observe that formal changes of coordinates
in $\bC^N$ and $\bC^{N'}$ give a
polynomial change of coordinates in
$J^k(\bC^N,\bC^{N'})_{(0,0)}$.

The
reader is referred e.g. to [GG] or [BER4,
Chapter XII] for further discussion of
these notions. Our main result is the following, which in
particular implies Theorem 2.1.1 above and Theorem 1 in the
introduction (we leave the details of the proofs of these implications
to the reader), and from which several other
theorems will be deduced below.  

\proclaim{Theorem 2.1.5} Let $M$ and $M'$
be formal generic submanifolds through
$0\in\bC^N$ and $0\in\bC^{N'}$ of
codimension $d$ and $d'$, respectively,
such that
$M$ is of finite type at $0$ and $M'$ is
$\ell_0$-nondegenerate at $0$ for some
integer
$\ell_0$. Assume that
$n\geq n'$, where
$n=N-d$ and $n'=N'-d'$.
Then there exist an integer $k_1$ with
$1 < k_1\leq d+1$, a polynomial
$P$ on $J^1(\bC^{n'},\bC^{n'})_{(0,0)}$ and for
each $\tilde \jmath=(j_1,\ldots, j_{n'})$, with
$1\leq j_1<\ldots <j_{n'}\leq n$, a linear
surjective mapping 
$$
\pi_{\tilde \jmath}\:
J^1(\bC^N,\bC^{N'})_{(0,0)}\to
J^1(\bC^{n'},\bC^{n'})_{(0,0)}
$$
and a formal power series in
$Z=(Z_1,\ldots Z_N)$
$$
\Phi^{\tilde
\jmath}(Z,\Lambda)\sim\sum_{|\alpha|>0}\frac{
c^{\tilde \jmath}_\alpha(\Lambda)}{(P(\pi_{\tilde
\jmath}(j_0^{k_1\ell_0,1}(\Lambda))))^
{l^{\tilde \jmath}_\alpha}} Z^\alpha, \tag2.1.6
$$
where $c^{\tilde \jmath}_\alpha(\Lambda)$ are $\bC^{N'}$
valued polynomials in
$J^{k_1\ell_0}(\bC^N,\bC^{N'})_{(0,0)}$ and
$l^{\tilde \jmath}_\alpha$ nonnegative integers,
satisfying the following. For every formal
mapping 
$H\:(M,0)\to(M',0)$, which is CR submersive i.e.
$$
d\Cal H(\Cal D_M(0))=\Cal
D_{M'}(0),\tag2.1.7
$$
there exists $\tilde \jmath$ as above such that
$P(\pi_{\tilde \jmath}(j_0^{1}(H)))\neq 0$ and 
$$
\aligned
H(Z)&\sim \Phi^{\tilde
\jmath}(Z,j^{k_1\ell_0}_0(H)),\ \text{\rm if
$k_1$ is even,}\\ H(Z)&\sim \Phi^{\tilde
\jmath}\left(Z,\overline{j^{k_1\ell_0}_0(H)}
\right),\ \text{\rm if
$k_1$ is odd}.
\endaligned\tag2.1.8
$$
In addition, if $M$ and $M'$ are
real-analytic, then for every $\tilde \jmath$
as above and $\Lambda_0$ with 
$P(\pi_{\tilde
\jmath}(j_0^{k_1\ell_0,1}(\Lambda_0))\neq 0$
the series \thetag{2.1.6} converges
uniformly for $(Z,\Lambda)$ near
$(0,\Lambda_0)$ in
$\bC^N\times
J^{k_1\ell_0}(\bC^N,\bC^{N'})_{(0,0)}$.
\endproclaim

In
what follows, we shall denote by
$\hat \Cal F(M,M')=\hat \Cal F(M,0;M',0)$ the set of formal
mappings $(M,0)\to(M',0)$ which are CR submersive. (For brevity,
we here suppress the dependence on the base points of $M$ and $M'$, which in this
section are assumed to the  origin in the respective
spaces.) When $M$ and $M'$
are 
real-analytic, then we also denote by
$\Cal F(M,M'):= \Cal F(M,0;M',0)$ those formal mappings
in
$\hat\Cal F(M,M')$ that are convergent,
and hence define holomorphic mappings 
which map a neighborhood of
$0$ in
$M$ into a neighborhood of $0$ in $M'$. Thus, in this notation, Theorem 3 in
the introduction gives sufficient conditions on $M$ and $M'$ so that
$\Cal F(M,M')=\hat\Cal F(M,M')$.

The following result, which will be proved
in \S4.3, is based on Theorem 2.1.5. 

\proclaim{Theorem 2.1.9} Let $M$ and $M'$
be formal generic submanifolds through
$0\in\bC^N$ and $0\in\bC^{N'}$,
respectively, such that
$M$ is of finite type at $0$ and $M'$ is
$\ell_0$-nondegenerate at $0$ for some
integer
$\ell_0$.  Then there exist an integer
$k_1$, depending only on $M$, with $1<k_1\leq d+1$ where $d$
denotes the codimension of $M$, and two real
algebraic subvarieties
$A, B\subset
J^{k_1\ell_0}(\bC^{N},\bC^{N'})_{(0,0)}$
such that the image of the mapping
$$
j^{k_1\ell_0}_0\:\hat\Cal F(M,M')\to
J^{k_1\ell_0}(\bC^{N},\bC^{N'})_{(0,0)}
\tag2.1.10
$$
coincides with $A\setminus B$. In
addition, the image $A\setminus B$ is
totally real at each nonsingular point.  
\endproclaim
\remark{Remark $2.1.11$} If $M$ and $M'$
are real-analytic, then, in view of Theorem
3, the conclusion of Theorem 2.1.9
also holds for convergent maps, i.e. with
$\hat\Cal F(M,M')$ replaced by $\Cal
F(M,M')$ in \thetag{2.1.10}. In this case, \thetag{2.1.10} is a
homeomorphism onto its image. This is the content of
Theorem 4, which will be proved in \S4.3.
\endremark\medskip

Let us consider the case $N=N'$ and
$M=M'$, where $M$ is of finite type at
$0$. By Corollary 1.27, the set
$\hat\Cal F(M,M)$ consists of formal
mappings $H\:(M,0)\to (M,0)$ which are
invertible. Thus, $\hat\Cal F(M,M)$
is a group under composition. 

For $k\geq 1$, we denote by
$G^k(\bC^N)_0$ the group (under
composition) of invertible jets in
$J^k(\bC^N,\bC^N)_{(0,0)}$, which is a
complex Lie group. It follows from the
above that, for any $k\geq 1$, the image of
$\hat\Cal F(M,M)$ under $j^k_0$ is
contained in $G^k(\bC^N)_0\subset
J^k(\bC^N,\bC^N)_{(0,0)}$.

\proclaim{Theorem 2.1.12} Let $M$ be a
formal generic submanifold through $0$ in
$\bC^N$ which is $\ell_0$-nondegenerate
and of finite type at $0$. Then there
exists an integer $k_1$ with $1<k_1\leq
d+1$, where $d$ denotes the codimension of
$M$, such that
$$
j^{k_1\ell_0}_0\:\hat\Cal F(M,M)\to
G^{k_1\ell_0}(\bC^N)_0\tag2.1.13
$$
is an injective group homomorphism and its
image is a totally real algebraic Lie
subgroup of $G^{k_1\ell_0}(\bC^N)_0$.
\endproclaim

If $M$ is a real-analytic generic
submanifold through $0$ in $\bC^N$, then, as mentioned in the
introduction, 
$\Cal F(M,M)$ (which in view of Theorem
3 coincides with $\hat\Cal F(M,M)$)
is usually called the stability
group of $M$ at $0$, and is denoted by
$\Aut(M,0)$. The group $\Aut(M,0)$ has a
natural (inductive limit) topology
corresponding to uniform convergence on
compact neighborhoods of $0$. That is, a
sequence $\{H_j\}
\subset \Aut(M,0)$ converges to $H \in
\Aut(M,0)$ if there is a compact
neighborhood  of $0$ to which all the
$H_j$ extend and on which the
$H_j$ converge uniformly to $H$.  

\proclaim{Theorem 2.1.14} Let $M$ be a
real-analytic generic submanifold through
$0$ in
$\bC^N$ which is $\ell_0$-nondegenerate
and of finite type at $0$. Then there
exists an integer $k_1$ with $1<k_1\leq
d+1$, where $d$ is the codimension of $M$,
such that
$$
j^{k_1\ell_0}_0\:\Aut(M,0)\to
G^{k_1\ell_0}(\bC^N)_0\tag2.1.15
$$
is a continuous injective group
homomorphism and its image is a totally
real algebraic Lie subgroup of
$G^{k_1\ell_0}(\bC^N)_0$. Moreover,
\thetag{2.1.15} is a homeomorphism onto
the image $j^{k_1\ell_0}_0(\Aut(M,0))$. 
\endproclaim

For the proofs of the results above, we shall need several
tools which will be presented in \S3 below. However, we first discuss
briefly the necessity of the conditions imposed on $M$ and $M'$ in the
results above.

\subhead 2.2. Generic necessity of finite type and 
finite nondegeneracy in the real-analytic
case\endsubhead In the theorems given in \S2.1, a
standing  assumption is that $M$ is of finite type at
$0$ and  that $M'$ is finitely nondegenerate at $0$.
In this section, we shall discuss to what extent these
conditions are necessary for the results. More
precisely, we shall discuss failure of the conclusion
in Theorem 2.1.1 (which is a consequence of the main
result, Theorem 2.1.5). 

The notion of finite
nondegeneracy at a point $p$ in a real-analytic, 
generic submanifold
$M$ is intimately related to that of holomorphic 
nondegeneracy as defined in the introduction. A 
connected, real-analytic, generic submanifold
$M\subset \bC^N$ of codimension $d$ is holomorphically
nondegenerate (at some point or, equivalently, at all
points) if and only there exists $\ell(M)$, $0\leq
\ell(M)\leq N-d$, such that $M$ is
$\ell(M)$-nondegenerate outside a proper real-analytic
subvariety of $M$ (see e.g. [BER1] or [BER4, Chapter
XI]). Also, it is easy to see that the set  of points
at which a real-analytic, generic submanifold is not
of finite type is a real-analytic subvariety of $M$
(see also [BER4, \S 1.5]). Thus, a connected,
real-analytic, generic submanifold $M\subset \bC^N$ of
codimension $d$ is either (a) $\ell$-nondegenerate,
for some $\ell$ with $0\leq \ell\leq N-d$, and of
finite type outside a proper  real-analytic subvariety
of $M$, (b) holomorphically degenerate, or (c) of
infinite type at every point (but (b) and (c) are not
mutually exclusive). 

For a formal generic submanifold $M\subset\bC^N$, the
notion of (formal) holomorphic nondegeneracy can be
defined as follows. We say that a formal $(1,0)$
vector field is {\it (formally) holomorphic} is its
coefficients are independent of $\zeta$. The formal
generic submanifold $M$ is said to be {\it
holomorphically nondegenerate} at $0$ if there are no
nontrivial (formal) holomorphic vector fields tangent
to $M$. If $M$ is a real-analytic generic submanifold,
then it is (formally) holomorphically nondegenerate at
0 as a formal submanifold if and only if it is
holomorphically nondegenerate at 0 as a real-analytic
one (i.e.\ in the sense defined in the introduction).
Moreover, if $M$ is a smooth generic submanifold,
then it is holomorphically nondegenerate at 0 (as a
formal submanifold) if and only there exists a
sequence of points $p_j\in M$ tending to $0$ such that
$M$ is finitely nondegenerate at each $p_j$. The
reader is referred to [BER4, Chapter XI] for these
results. 

The following result shows necessity of the
hypotheses in Theorem 2.1.1.

\proclaim {Theorem 2.2.1} Let $M\subset\bC^N$ be a 
formal generic submanifold. Suppose either of the
following hold. 
\roster
\item "(i)" $M$ is holomorphically
degenerate at $0$.
\item "(ii)" $M$ is weighted
homogeneous, i.e. defined by the 
vanishing of weighted homogeneous polynomials, and of
infinite type at every point (as a real-analytic
submanifold).
\endroster
Then for any integer $K >
0$ there exist local formal invertible mappings
$$H^1,H^2\:(\bC^N,0)\to (\bC^N,0)$$
mapping $M$ into itself such that
$$
\frac{\partial^{|\alpha|}H^1}{\partial Z^\alpha}(0)=
\frac{\partial^{|\alpha|}H^2}{\partial
Z^\alpha}(0),\quad \forall  |\alpha|\leq K,\tag2.2.2
$$ 
but $H^1 \not\equiv H^2$. If $M$ is real-analytic,
then $H_1$ and $H_2$ can be chosen to be biholomorphic
near $0$.
\endproclaim

In the real-analytic case, Theorem 2.2.1 was proved in
[BER2]. The proof of Theorem 2.2.1 in the general case
is similar to that in the real-analytic case and the
modifications are left to the reader. 

\heading \S3. Tools for the proofs\endheading

\subhead \S3.1. The Segre mappings\endsubhead We keep the notation
introduced in the previous 
sections; e.g.\ $M$ is a formal generic
submanifold of codimension
$d$ defined by the formal power series $\rho=(\rho_1,\ldots,
\rho_d)$. Recall that $Z=(z,w)$, with $z=(z_1,\ldots,
z_n)$ and $w=(w_1,\ldots, w_d)$, and $\zeta=(\chi,\tau)$, with $ \chi=(\chi_1,\ldots,
\chi_n)$ and $\tau= (\tau_1,\ldots \tau_d)$, are normal coordinates
for $M$ at $0$, so that $M$ is defined by
$$
w_j-Q_j(z,\chi,\tau),\quad j=1,\ldots,d,\tag3.1.1
$$
where the
$Q_j\in\bC[[z,\chi,\tau]]$ satisfy
$$
Q_j(0,\chi,\tau)\sim
Q_j(z,0,\tau)\sim\tau_j.\tag3.1.2
$$
As mentioned in \S1, $M$ is also defined by 
$$
\tau_j-\bar Q_j(\chi,z,w),\quad
j=1,\ldots,d.\tag3.1.3
$$
Consider, for each integer $k\geq1$, the formal
mapping $v^k\:(\bC^{kn},0)\to(\bC^{N},0)$
defined as follows. For $k=2j$, 
$$
\multline
v^{2j}(z,\chi^1,\ldots,z^{j-1},\chi^j)=\bigg(z,
Q\big(z,\chi^1,\bar
Q\big(\chi^1,z^{1},Q\big(z^{1},
\chi^{2},\ldots,\\\bar
Q\big(\chi^{j-1},z^{j-1},Q\big(z^{j-1},
\chi^j,0\big)\big)\ldots\big)
\big)\big)\bigg),
\endmultline\tag3.1.4
$$  
and, for $k=2j+1$, 
$$
\multline
v^{2j+1}(z,\chi^1,\ldots,z^{j-1},\chi^j,
z^j)=\bigg(z,Q\big(z,\chi^1,\bar
Q\big(\chi^1,z^{1},Q\big(z^{1},\chi^{2},
\ldots,\\ Q\big(z^{j-1},\chi^j,\bar
Q\big(\chi^j,z^j,
0\big)\big)\ldots\big)\big)\big)\bigg).
\endmultline\tag3.1.5
$$
(In \thetag{3.1.4} and \thetag{3.1.5}, $z^r$ denotes
$(z^r_1,\ldots,z^r_n)$ and similarly
$\chi^r$ denotes $(\chi^r_1,\ldots,
\chi^r_n)$ for $r=1,\ldots,
j$.) For $k=0$, we set $v^0=(0,0)$. We shall
refer to the mapping $v^k$ as the
$k$th {\it Segre mapping} of $M$. 

\proclaim{Proposition 3.1.6} For any defining series $\rho\in\bC[[Z,\zeta]]^d$ of
$M$ and any $k\geq 0$,
$$
\rho(v^{k+1}(z,\chi^1, z^1,\ldots),\bar
v^k(\chi^1,z^1,\ldots))\sim 0.\tag3.1.7
$$
\endproclaim

\demo{Proof} For simplicity, we only
consider the case where $k=2j$. It follows
from
\thetag{3.1.4} and
\thetag{3.1.5} that
$$
v^{2j+1}(z,\chi^1,\ldots,z^{j-1},\chi^j,
z^j)\sim\left (z,Q(z,\bar
v^{2j}(\chi^1,\ldots,z^{j-1},\chi^j,
z^j)\right).\tag3.1.8
$$
It suffices to show \thetag{3.1.7} for the defining series 
given by 
\thetag{3.1.1}, for which \thetag{3.1.7}
is an immediate consequence of
\thetag{3.1.8}. This completes the proof.\qed 
\enddemo

The following characterization of finite type will be 
important.  If $v$ is a formal mapping $(\bC^m,0)\to
(\bC^l,0)$, then we shall write 
$\Rk (v)$ to denote the rank of the matrix $(\partial
v_i/\partial x_j)$,
$i=1,\ldots, l$, $j=1,\ldots m$, in $\Bbb A^l$, where
$\Bbb A$ denotes the field of
fractions of $\bC[[x_1,\ldots, x_m]]$. We shall also the
notation $\rk(\partial
v_i/\partial x_j)$ for this rank. 

\proclaim{Theorem 3.1.9} Let $M$ be a
formal generic submanifold of $\bC^N$
through $0$. Then, the following are
equivalent:
\roster
\item"(i)" $M$ is of finite type at $0$;
\item"(ii)" There exists $k_1\leq d+1$ such
that the rank
$\Rk(v^k)$ is
$N$ for
$k\geq k_1$.
\endroster
\endproclaim

For the proof of Theorem 3.1.9, we shall need special coordinates for
a formal generic submanifold. These will be presented in \S3.2
below. The proof of Theorem 3.1.9 will be given in \S3.3.

\subhead \S3.2. Formal canonical
coordinates\endsubhead 
Recall that $\Cal D_M$ denotes
the $\bC[[Z,\zeta]]$-module generated by all the formal
$(1,0)$ and $(0,1)$ vector fields tangent to the
formal generic submanifold $M$ of $\bC^N$. We define
the integers
$m_1,\ldots,m_h$, also called the  {\it
H\"ormander numbers} of
$M$ at $0$, as follows. The number $m_1$ is the
smallest integer for which there exists a
commutator $C$ of vector fields in $\Cal D_M$ of
length
$m_1$\footnote{The length of a
commutator is the number of vector fields in
$\Cal D_M$ appearing; e.g., the commutator
$[X,[Y,Z]]$, with $X,Y,Z\in\Cal D_M$, has length
$3$.}  such that
$C(0)$ is not in the span of $\Cal
D_M(0)\subset T'_{0}\bC^{2N}$. We define the
subspace
$E_1\subset T'_{0}\bC^{2N}$ to be the linear
span of $\Cal D_M(0)$ and the values at
$0$ of all commutators of vector fields in
$\Cal D_M$ of length $m_1$. We define $l_1$ to
be
$$ 
l_1=\dim E_1-2(N-d).
\tag3.2.1
$$ 
We define inductively the numbers
$m_1<m_2<\ldots<m_h$ and subspaces
$E_1\subset E_2\subset\ldots\subset
E_h=T'_{0}\bC^{2N}$ as follows. The number
$m_{k+1}$ is the smallest integer
for which there exists a commutator $C$ of
vector fields in $\Cal D_M$ of length $m_{k+1}$
such that $C(0)\not\in E_k$. The subspace
$E_{k+1}$ is then defined as the span of $E_k$
and the values at
$0$ of all commutators of vector fields in
$\Cal D_M$ of length
$m_{k+1}$. We define
$$ 
l_{k+1}=\dim E_{k+1}-\dim E_k=\dim
E_{k+1}-2(N-d)-\sum_{i=1}^k l_i.\tag 3.2.2
$$ 
It is clear that this process terminates
after a finite number of steps.
We shall call the number $l_j$ the {\it
multiplicity} of the H\"ormander number
$m_j$. It is also convenient to use the notation
$\mu_1,\ldots,\mu_r$ for the H\"ormander numbers
repeated according to their multiplicities, so
that $r=\sum_{j=1}^h l_j$. 

The following theorem will be used in the
proof of Theorem 3.1.9. 

\proclaim{Theorem 3.2.3} Let $M$ be
a formal generic submanifold of
$\bC^N$ of codimension
$d$ through $0$. Let $2\le
\mu_1\leq\ldots\leq\mu_r$ be the H\"ormander
numbers of $M$ at $0$ repeated according to their
multiplicities. There exists a formal change of
coordinates
$Z=Z(z,w',w'')$, with
$z=(z_1,\ldots, z_n)$, $w'=(w'_1,\ldots,w'_r)$,
$w''=(w''_1,\ldots,w''_{d-r})$, $N=n+d$,
satisfying the following. The defining series $\rho$ 
of $M$, after the formal change of coordinates
$(Z,\zeta)=(Z(z,w',w''),\bar Z(\chi,\tau',\tau''))$, satisfies 
$$
\rho(z,w',w'',\chi,\tau',\tau'')=a(z,w',w'',\chi,\tau',\tau'')\pmatrix
w'-Q'(z,\chi,\tau',\tau'')\\w''-Q''(z,\chi,\tau',\tau'')\endpmatrix\tag3.2.4
$$
where $a(z,w',w'',\chi,\tau',\tau'')$ is a $d\times d$ matrix of
formal power series which is 
invertible at $0$, and $Q'(z,\chi,\tau',\tau'')$
and
$Q''(z,\chi,\tau',\tau'')$ are of the form
$$
\aligned
Q'_k(z,\chi,\tau',\tau'')\sim&\, \tau'_k +p_k(z,\chi,
\tau'_1,\ldots,\tau'_{k-1})+A_k(z,\chi
,\tau',\tau'')\tau''+\\&\hskip 2.5in
R_k(z,\chi,\tau')\\ Q''(z,\chi,\tau',\tau'')\sim&\,
\tau'' + B(z,\chi,\tau',\tau'')\tau'',
\endaligned\tag 3.2.5
$$ 
where $k=1,\ldots, r$.
Here, $p_k(z,\chi, \tau'_1,\ldots,\tau'_{k-1})$ is
a weighted homogeneous polynomial of
degree $\mu_k$, where
$z$ and
$\chi$ have weight one and $\tau'_j$ has weight $\mu_j$
for $j=1,\ldots r$, $R_k(z,\chi,\tau')$ is a 
formal power series which is
$O(\mu_k+1)$ (i.e. involving only terms which are
weighted homogeneous of degree at least $\mu_k+1$),
$A_k(z,\chi,\tau',\tau'')$ and
$B(z,\chi,\tau',\tau'')$ are matrices of formal
power series without constant terms. Moreover, we
have
$$
\aligned
&Q'(z,0,\tau',\tau'')\sim
Q'(0,\chi,\tau',\tau'')\sim
\tau',\\ &Q''(z,0,\tau',\tau'')\sim
Q''(0,\chi,\tau',\tau'')\sim\tau''.\endaligned\tag3.2.6
$$ 
\endproclaim

The proof of Theorem 3.2.3 can be extracted from 
the proof of [BER4, Theorem 4.5.1]. The reader
should observe that $M$ is of finite type at $0$
if and only if $r=d$. In this case, there are no 
$w''$ variables in Theorem 3.2.3, i.e. $w=w'$ and
$\tau=\tau'$ in \thetag{3.2.4--6}.

\subhead \S3.3. Proof of Theorem 3.1.9\endsubhead  In view of Theorem
3.2.3, we may assume
that we have formal coordinates
$Z=(z,w',w'')$,
$\zeta=(\chi,\tau',\tau'')$, as described in
Theorem 3.2.3, such that $M$ is defined by $w'-Q'(z,\chi,\tau',\tau'')$ and
$w''-Q''(z,\chi,\tau',\tau'')$, where $Q'$ and
$Q''$ satisfy \thetag{3.2.5} and \thetag{3.2.6}. Let
us write 
$$
Q'(z,\chi,\tau',\tau'')\sim
\tau'+p(z,\chi,\tau')+R(z,\chi,\tau,\tau''),
\tag3.3.1
$$
where 
$$
p(z,\chi,\tau')=(p_1(z,\chi),\ldots,
p_r(z,\chi,\tau'_1,\ldots,\tau'_{r-1}))\tag3.3.2
$$ 
are
weighted homogeneous polynomials as in Theorem
3.2.3 and
$R=(R_1,\ldots R_r)$ are the remainder terms of
higher (weighted) homogeneity. Consider the homogeneous
generic submanifold $M^0$ of $\bC^N$ given by
$$
w'=\bar w'+p(z,\bar z,\bar w'),\quad
w''=\bar w''.\tag3.3.3
$$
Observe that $M^0$ has the same H\"ormander
numbers as $M$ (with multiplicity). For each
fixed $k$, we denote, for simplicity of notation,
the variables in the space
$\bC^{kn}$, where the $k$th Segre mappings are
defined, by
$(z,\xi)$, where $\xi\in\bC^{(k-1)n}$.  We denote
by
$v^k(z,\xi)$ the $k$th Segre mapping
of $M$ at $0$ as defined by \thetag{3.1.4} and
\thetag{3.1.5}, and by
$v^k_0(z,\xi)$ the
$k$th Segre mapping of the formal generic
submanifold associated to $M^0$ at $0$. The $j$th
component of the mapping $v^k_0$ is a 
homogeneous polynomial (in the usual sense; i.e.\ all
components of $z$ and $\xi$ have weight one) of degree 
$\mu_j$, where
$\mu_j$ denotes the $j$th H\"ormander number
(with multiplicity) of
$M^0$ (or $M$) at $0$. The
$z$ component of the mappings $v^k(z,\xi)$ and
$v^k_0(z,\xi)$ coincide, and are equal to $z$.
Moreover, the
$w''$ components also coincide, and are equal to
$0$. Let
us separate the $z$,
$w'$, and $w''$ components of the mappings $v^k$
and $v^k_0$ and write
$$
v^k(z,\xi)=(z,g^k(z,\xi),0),\quad
v_0^k(z,\xi)=(z,g_0^k(z,\xi),0),\tag3.3.4
$$
where $(g^k_0)_j(z,\xi)$ is a homogeneous
polynomial of degree $\mu_j$. We have, for
$j=1,\ldots, r$,
$$
g^k_j(z,\xi)\sim(g^k_0)_j(z,\xi)+O(\mu_j+1),
\tag3.3.5
$$
where $O(\nu+1)$ denotes a power series consisting
only of terms of degrees higher than $\nu$. 

It follows from [BER1, Proposition 2.4.1] (see
also [BER4, Proposition 10.5.27]) that there
exists $k_1$, with $1\leq 
k_1\leq d+1$, such that
$\Rk (v^k_0)=n+r$ for $k\geq k_1$. Observe
that  the determinant of an $m\times m$
matrix
$A$ of power series, where the $j$th row
of $A$ is of the form $f_j=f^0_j+O(d_j+1)$
for some homogeneous polynomial $f^0_j$ of
degree $d_j$, is of the form
$$
\det\, A=\det\, A^0+O(d_1+\ldots+d_m+1),\tag3.3.6
$$
where $A^0$ is the matrix with rows $f^0_1,\ldots
f^0_m$. It follows from this observation and
\thetag{3.3.5} that
$\Rk(v^k)\geq \Rk (v^k_0)$, and hence 
$\Rk (v^k)\geq n+r$ for $k\geq k_1$. On
the other hand, by the form \thetag{3.3.4}
of $v^k$, we have
$\Rk (v^k)\leq n+r$ for any $k\geq 1$. Thus, 
$$
\Rk (v^k)=n+r,\quad \forall k\geq k_1,\tag3.3.7
$$
The equivalence of (i) and (ii) of Theorem 3.1.9 follows from the fact
that $M$ is of finite type at $0$ if and only if $r=d$, i.e. if and
only if $\Rk (v^k)=n+d=N$, for $k\geq k_1$. \qed

\subhead \S3.4. Basic identity for formal
mappings\endsubhead
An important tool, in combination with the Segre mappings, in the
proofs of the 
theorems in \S2.1 will be the basic
identity which we shall present in this
section. We keep the notation
established in the previous sections. In what follows, $M$ and $M'$ 
denote fixed formal generic submanifolds
of codimension $d$ and $d'$ through the
origin of $\bC^N$ and $\bC^{N'}$,
respectively. 

We
choose normal coordinates $Z=(z,w)$,
$\zeta=(\chi,\tau)$ for $M$ as in \S1,
and, similarly, normal coordinates
$Z'=(z',w')$,
$\zeta'=(\chi',\tau')$ for $M'$. There is
an associated coordinate system on
$J^l(\bC^N,\bC^{N'})_{(0,0)}\cong
\bC^{K(l)}$, where $K(l)$ denotes the
dimension of this jet space. We shall
use a scaled coordinate system whose
coordinates we shall denote  by
$$
\Lambda=\left(\lambda_{z^\alpha
w^\beta},\mu_{z^\gamma
w^\delta}\right)_{1\leq|\alpha|+|\beta|,
|\gamma|+|\delta|\leq l}.\tag3.4.1
$$
For a formal mapping $H\:(\bC^N,0) \to
(\bC^{N'},0)$, 
write
$H(z,w)=(F(z,w),G(z,w))$, where
$F=(F_1,\ldots,F_{n'})$ and $G=(G_1,\ldots
G_{d'})$. In the scaled coordinates
\thetag{3.4.1} we have, for any
non-negative integer $l$,
$$
j_0^l(H)=\left(\lambda_{z^\alpha
w^\beta},\mu_{z^\gamma
w^\delta}\right)_{1\leq|\alpha|+|\beta|,|\gamma|+|\delta|\leq
l},\tag3.4.2
$$
where 
$$
\lambda_
{z^\alpha
w^\beta}=\frac
{\partial^{|\alpha|+|\beta|}F}{\partial
z^\alpha\partial w^\beta}(0,0),\quad
\mu_ {z^\gamma
w^\delta}=\frac
{\partial^{|\gamma|+|\delta|}G}{\partial
z^\gamma\partial w^\delta}(0,0).\tag3.4.3
$$ 
For each fixed
$l$, we shall split, and reorder, the
variables
$\Lambda$ in \thetag{3.4.1} as follows
$$
\Lambda=(\Lambda',\Lambda'')\tag3.4.4
$$
where 
$$
\Lambda''=\left(\mu_{z^\gamma}
\right)_{1\leq|\gamma|\leq l}\tag3.4.5
$$
and the components of $\Lambda'$ are the 
remaining variables in \thetag{3.4.1}. We
shall denote the number of components of
$\Lambda'$ by $K'=K'(l)$ and that of
$\Lambda''$ by $K''=K''(l)$, so that
$J^l(\bC^N,\bC^{N'})_{(0,0)}\cong
\bC^{K'}\times \bC^{K''}$.  We are now in
a position to state the basic identity.

\proclaim{Theorem 3.4.6} Let $M$ and $M'$
be formal generic submanifolds through
the origin in $\bC^N$ and $\bC^{N'}$,
respectively.  Assume that $M'$ is
$\ell_0$-nondegenerate at $0$, and that
$n\geq n'$, where $n=N-d$, $n'=N'-d'$,
$d=\codim M$, and $d'=\codim M'$. Then for
each
$\tilde \jmath=(j_1,\ldots, j_{n'})$, with
$1\leq j_1<\ldots <j_{n'}\leq n$, 
and for every $\alpha\in\Bbb
Z_+^{N}$, there exists a formal power
series mapping of the form 
$$
\Psi^{\tilde
\jmath}_\alpha(Z,\zeta,\zeta',\Lambda)
\sim\sum_
{\beta,\gamma,\delta,\kappa}\frac{
d_{\beta\gamma\delta\kappa}(\Lambda')}
{(\det\,(\lambda^l_{z_{j_p}})_{1\leq
l,p\leq n'} )^
{l_{\beta\gamma\delta\kappa}}}
Z^\beta\zeta^\gamma{\zeta'}^\delta
{\Lambda''}^\kappa,
\tag3.4.7
$$
where  
$\Lambda=(\Lambda',\Lambda'')\in
\bC^{K'}\times \bC^{K''}$ with
$K'=K'(\ell_0+|\alpha|)$ and
$K''=K''(\ell_0+|\alpha|)$,
$d_{\beta\gamma\delta\mu}(\Lambda')$
are
$\bC^{N'}$ valued polynomials in
$\bC^{K'}$,
and
$l_{\beta\gamma\delta\mu}$ nonnegative
integers, satisfying the following. For
every formal mapping $H\in\hat\Cal
F(M,M')$ there exists $\tilde \jmath$ such that 
$$
\det\,\left(\frac{\partial
F_l}{\partial z_{j_p}}(0)\right)_{1\leq l,
p\leq n'}\neq 0\tag3.4.8
$$ 
and, for all $\alpha\in\Bbb Z^N_+$,
$$
\partial^\alpha
H(Z)-\Psi^{\tilde \jmath}_{\alpha}\left
(Z,\zeta,\bar H(\zeta),(\partial^\beta\bar
H(\zeta))_{1\leq |\beta|\leq
\ell_0+|\alpha|}\right)\sim a(Z,\zeta)\rho(Z,\zeta),\tag3.4.9
$$
where $a(Z,\zeta)$ is a $d\times d$ matrix of formal power series and
$\rho=(\rho_1,\ldots,\rho_d)$ is a defining series for 
$M$. Moreover, \thetag{3.4.9} holds for
any $H\in\hat\Cal F(M,M')$ and any
$\tilde \jmath$ such that
\thetag{3.4.8} holds. 

If $M$ and $M'$ are real-analytic, then,
for any $\tilde \jmath$ as above, $\alpha\in\Bbb
Z^N_+$, and any 
$\Lambda'_0\in\bC^{K'(\ell_0+|\alpha|)}$
satisfying 
$$
\det\,\left({(\lambda_0)}^l_{z_{j_p}}\right)
_{1\leq l, p\leq n'}\neq 0,\tag3.4.10
$$ 
the series
$\Psi^{\tilde
\jmath}_\alpha(Z,\zeta,\zeta',\Lambda)$ given by
\thetag{3.4.7} defines a holomorphic
mapping near the point
$(Z,\zeta,\zeta',\Lambda',\Lambda'')=
(0,0,0,\Lambda_0,0)$.
\endproclaim
\remark{Remark $3.4.11$} The reader should
observe that in substituting the
formal mapping $(\partial^\beta\bar
H(\zeta))_{1\leq |\beta|\leq
\ell_0+|\alpha|}$ for $\Lambda$ in
\thetag{3.4.7}, we replace $\Lambda''$ by 
$(\partial_\chi^\beta\bar
G(\chi,\tau))_{1\leq |\beta|\leq
\ell_0+|\alpha|}$ and $\Lambda'$ by the
remaining derivatives. This substitution
of formal power series makes sense since,
as remarked in \S1, $G(z,0)\sim 0$, and
the dependence on $\Lambda'$ is rational.
In what follows, we shall, for fixed $l$,
decompose, and reorder the components of 
$(\partial^\beta
H(Z))=(\partial^\beta
H(Z))_{1\leq |\beta|\leq l}$ as
$((\partial^\beta
H(Z))',(\partial^\beta
H(Z))'')$, where 
$(\partial^\beta
H(Z))''=(\partial_z^\beta
G(z,w))_{1\leq |\beta|\leq l}$ and
$(\partial^\beta H(Z))'$  denotes the
remaining derivatives. 
\endremark\medskip

\demo{Proof of Theorem $3.4.6$}  Recall that, in the chosen
normal coordinates,
$H$ maps
$(M,0)\to (M',0)$ if and only if
\thetag{1.14} holds. The mapping $H$
belongs to $\hat\Cal F(M,M')$ if it also
satisfies \thetag{1.19}, which, as noted
in the proof of Corollary 1.27, is
equivalent to $\partial F/\partial z(0,0)$
having rank $n'$. Thus, there exists
$\tilde \jmath$ as in the statement of the
theorem such that \thetag{3.4.8} holds. 

Let 
$\tilde \jmath$, as in the theorem, be given.
We shall consider only those mappings
$H\:(M,0)\to (M',0)$ for which
\thetag{3.4.8} holds.  
After renumbering the variables
if necessary, we may assume that
$\tilde \jmath=(1,2,\ldots, n')$.

We take as a basis for the $(0,1)$ vector
fields tangent to $M$ the following
$$
L_j=\frac{\partial}{\partial
\chi_j}+\sum_{k=1}^d \bar Q_{k,\chi_j}(\chi,
z,w)\frac{\partial}{\partial
\tau_k}, \ \ j = 1,\ldots, n.\tag3.4.12
$$
 In what follows, we
shall only use the vector fields
$L_1,\ldots, L_{n'}$. We shall also need
the following vector fields tangent to
$M$ and given by
$$
\aligned \tilde
L_j&=\frac{\partial}{\partial
z_j}+\sum_{k=1}^d  Q_{k,z_j}(z,
\chi,\tau)\frac{\partial}{\partial
w_k}, \quad\quad j=1,...,n,\\ 
T_j&=\frac{\partial }{\partial w_j}+\sum_{k=1}^d
\bar Q_{k,w_j}(\chi,
z,w)
\frac{\partial}{\partial
\tau_k},\quad\quad j=1,...,d,\\
V_j&=\tilde L_j-\sum_{k=1}^d
Q_{k,z_j}(z,\chi,
\tau)T_k,  \quad\quad j=1,...,n.
\endaligned\tag3.4.13
$$
Note that the $\tilde L_j$ form a basis
for the $(1,0)$ vector fields tangent to
$M$ (modulo those whose coefficients are in the ideal
generated by a set of defining power series of $M$).  

After applying
the $L_j$, for $j=1,\ldots, n'$, to the
second set of equations in \thetag{1.14}
$|\a|$ times, and applying Cramer's rule
after each application, we obtain, for
$w=Q(z,\chi,\tau)$ or $\tau=\bar
Q(\chi,z,w)$, any
multi-index
$\alpha$, and $l=1,\ldots, d'$,
$$
\multline
\bar Q_{l,{\chi'}^\a}'(\bar F(\chi,\tau), 
F(z,w),G(z,w))=\\
\sum_{1\leq |\beta|\le |\a|} 
(L^\beta\bar G_l (\chi,\tau))
P_{\a,\beta}\left((L^\gamma
\bar F(\chi,\tau))_{1\leq |\gamma|\le
|\a|}\right)/\Delta^{2|\a|-1},\endmultline
\tag 3.4.14
$$
where $\Delta =
\Delta(z,w,\chi,\tau)=
\det[L_j\bar F_k (\chi,\tau)]_{1\leq
j,k\leq n'}$, and
$P_{\a,\beta}$ are universal polynomials,
i.e.  independent of $M$,
$M'$, and
$H$. Note that for
any formal power series $h(\chi,\tau)$ and
any multi-index
$\beta$ we have
$L^\beta h(0)= \partial_{\chi^\beta}h(0)$.
Thus,
$\Delta(0)\not= 0$ by 
\thetag{3.4.8}. Also,
$L^\beta G(0) = 0$ by the normality of 
the coordinates, as is easily verified
from \thetag{1.14}. Hence the right hand
side of
\thetag{3.4.14} vanishes at the origin.

By the assumption that $M'$ is
$\ell_0$-nondegenerate at $0$, there exist
$n'$ multi-indices
$\a^1, \ldots, \a^{n'}$, with $1\le
|\a^j|\le
\ell_0$, and $n'$ integers
$l_1,\ldots, l_{n'}\in\{1,\ldots, d'\}$
such that
$\det\left[\bar
Q'_{l_j,{\chi'}^{\a^j}z'_k}(0)\right]\not= 0$. 
(See [BER1] or [BER4, Corollary 11.2.14].)
Hence, by the implicit function theorem,
there exists a unique $\bC^{n'}$ valued
formal power series $S(\chi',\tau',r)$
with $r=(r_1,\ldots, r_{n'})$,  so that
$S(0,0,0)=0$ and
$X=S(\chi',\tau',r)$ solves the system of
equations
$$
\bar
Q_{l_j,{\chi'}^{\a^j}}'
(\chi',X,Q'(X,\chi',\tau'))\sim r_j,
\quad j=1,\ldots,n'.  \tag 3.4.15
$$ 
For any
positive integer $k$, we shall introduce
the vector valued variables
$(a_\gamma)_{|\gamma|\leq k}$,
$(b_\beta)_{|\beta|\leq k}$, where
$\beta,\gamma\in\Bbb Z^{n'}_+$,
corresponding to $(L^\gamma
\bar F(\chi,\tau))_{|\gamma|\le
 k}$, $(L^\beta\bar G
(\chi,\tau))_{|\beta|\le
 k}$, respectively. Here
$a_\gamma=(a_\gamma^m)_{1\leq m\leq n'}$
and $b_\beta=(b^j_\beta)_{1\leq j\leq
d'}$. We write
$(a^m_{k})_{1\leq k,m\leq n'}$ for
$(a^m_\gamma)_{1\leq m\leq n',\,
|\gamma|=1}$. We define
$$
R_\alpha\left((a_\gamma)_{1\leq|\gamma|\leq|\alpha|},
(b^1_\beta)_{1\leq|\beta|\leq|\alpha|}\right)=\frac{\sum_{
1\leq|\beta|\le |\a|} (b^1_\beta)
P_{\a,\beta}\left((a_\gamma)_{1\leq|\gamma|\le
|\a|}\right)}{\left(\det (a^m_k)_{1\leq
k,m\leq n'}\right)^{2|\a|-1}}.\tag3.4.16
$$ 
Observe that $R_\alpha$ is a universal
rational function that vanishes when
$b^1_\beta=0$, $|\beta|\leq |\alpha|$, and
whose denominator is a power of 
$\det (a^m_k)_{1\leq k,m\leq n'}$. It
follows from the above that, for
$w=Q(z,\chi,\tau)$ or $\tau=\bar
Q(\chi,z,w)$, we have the identity
$$ 
F(z,w)=\Theta\left((L^\gamma
\bar F(\chi,\tau))_{|\gamma|\le
\ell_0},(L^\beta\bar G
(\chi,\tau))_{|\beta|\le
\ell_0}\right),\tag3.4.17
$$ 
where
$$
\aligned
\Theta((a_\gamma)_{|\gamma|\leq \ell_0 },&
(b_\beta)_{|\beta|\leq \ell_0})=\\&S
\left(a_0,b_0,\left(R_{\alpha^j}
\left((a_\gamma)_{1\leq |\gamma|\leq
|\alpha^j|}, (b^{l_j}_\beta)_{1\leq
|\beta|\leq
|\alpha^j|}\right)\right)_{1\leq j\leq
n'}\right).\endaligned\tag3.4.18
$$ 
Now, since
$F(z,w)$ is a power series in 
$(z,w)$ only, we have, for any multi-index
$\nu=(\nu^\prime,\nu^ {\prime\prime})$,
$$
V^{\nu^\prime}
T^{\nu^{\prime\prime}}F(z,w)=\frac{\partial^
{|\nu|}F}{\partial z^{\nu^\prime}\partial
w^{\nu^{\prime\prime}}} (z,w).\tag3.4.19
$$ 
By applying $V^{\nu^\prime}
T^{\nu^{\prime\prime}}$ to the identity
\thetag{3.4.17}, we obtain
$$
\aligned
\frac{\partial^ {|\nu|}F}{\partial
z^{\nu^\prime}\partial
w^{\nu^{\prime\prime}}}(z,w)=\\
\Theta_\nu\bigg((V^{\delta'}
T^{\delta''}
L^\gamma
\bar F(\chi,&\tau))_{|\delta|+|\gamma|\le
\ell_0+|\nu|},(V^{\kappa'}
T^{\kappa''}L^\beta\bar G
(\chi,\tau))_{|\kappa|+|\beta|\le
\ell_0+|\nu|}\bigg),\endaligned\tag3.4.20
$$ 
where we have used the notation
$\delta=(\delta',\delta'')$ and
$\kappa=(\kappa',\kappa'')$.   The identity
\thetag{3.4.20} holds when
$w=Q(z,\chi,\tau)$ or $\tau=\bar
Q(\chi,z,w)$. Observe that the power series $\Theta_\nu$
depends only on $\Theta$ and its derivatives. 

By substituting
\thetag{3.4.12} and \thetag{3.4.13} in
\thetag{3.4.20}, we obtain for
$w=Q(z,\chi,\tau)$ or $\tau=\bar
Q(\chi,z,w)$,
$$
\aligned
\frac{\partial^ {|\nu|}F}{\partial
z^{\nu^\prime}\partial
w^{\nu^{\prime\prime}}}(z,w)=\\
\Phi^1_\nu\bigl(
z,w,&\chi,\tau,\bar
F(\chi,\tau),\bar G(\chi,\tau),(\partial
^\alpha\bar H(\chi,\tau))', (\partial
^\alpha\bar
H(\chi,\tau))''\bigr).\endaligned\tag3.4.21
$$
where $(\partial
^\alpha\bar H(\chi,\tau))=(\partial
^\alpha\bar H(\chi,\tau))_{|\alpha|\leq
\ell_0+|\nu|}$ and we use the
notation $(\partial
^\alpha\bar H(\chi,\tau))=((\partial
^\alpha\bar H(\chi,\tau))',(\partial
^\alpha\bar H(\chi,\tau))'')$ as explained
in Remark 3.4.11. Observe that the power series $\Phi^1_\nu$
depends only on $M$ and $M'$, and not on the mapping $H$.
Using the notation
$$
\bar\Lambda=((\lambda_
{\chi^\alpha
\tau^\beta})_{1\leq |\alpha|+|\beta|\leq
\ell_0+|\nu|}, (\mu_{\chi^\gamma
\tau^\delta})_{1\leq|\gamma|+|\delta|\leq
\ell_0+|\nu|}),\tag3.4.22
$$
decomposing $\bar\Lambda$, and reordering
its components in an analogous
fashion as $\bar \Lambda=(\bar
\Lambda',\bar\Lambda'')$ with
$\bar\Lambda''=
(\mu_{\chi^\beta})_{1\leq|\beta|\leq
\ell_0+|\nu|}$,  it follows from
\thetag{3.4.18} and \thetag{3.4.20} that
the power series 
$$
\Phi^1_\nu\left(
z,w,\chi,\tau,\chi',\tau',\bar
\Lambda',\bar\Lambda''\right)\tag 3.4.23
$$ 
is of the form
$$
\Phi^1_\nu(z,w,\chi,\tau,\chi',\tau',\bar
\Lambda',\bar\Lambda'')
\sim\sum_
{\beta,\gamma,\delta,\kappa}\frac{
e_{\beta\gamma\delta\kappa}(\bar
\Lambda')} {\det(\lambda^j_{\chi_k})^
{l_{\beta\gamma\delta\kappa}}}
Z^\beta
\zeta^\gamma{\zeta'}^\beta 
\bar{\Lambda''}^\kappa
\tag3.4.24
$$
where  
$e_{\beta\gamma\delta\kappa}(\bar \Lambda)$
are
$\bC^{n'}$ valued polynomials and
$l_{\beta\gamma\delta\kappa}$ nonnegative
integers. We have used here the fact that
$$
R_{\alpha^j}
\left((a_\gamma)_{1\leq |\gamma|\leq
|\alpha^j|}, (b^{l_j}_\beta)_{1\leq
|\beta|\leq
|\alpha^j|}\right)=0,\tag3.4.25
$$
when
$(b^{l_j}_\beta)_{1\leq |\beta|\leq
|\alpha^j|}=0$. In view of
\thetag{3.4.21} and
\thetag{3.4.24}, we can take the first $n'$
components of $\Psi^{\tilde \jmath}_
\nu(Z,\zeta,\zeta', \Lambda)$, with
the fixed choice of
$\tilde \jmath$ above,  in the conclusion of
the theorem to be
$\Phi^1_\nu(z,w,\chi,\tau,\chi',\tau',
\Lambda',\Lambda'')$.  

To complete the construction
of
$\Psi^{\tilde \jmath}_\nu$, we need to find the
components corresponding to $G$ and its
derivatives. For this we substitute
\thetag{3.4.21} with
$\nu=0$ in the first set of equations in
\thetag{1.14}, and apply the
vector fields $V_j$ and $T_j$ to
the identity thus obtained, as above. We
obtain
$$
\aligned
\frac{\partial^ {|\nu|}G}{\partial
z^{\nu^\prime}\partial
w^{\nu^{\prime\prime}}}(z,w)=\\
\Phi^2_\nu\bigl(
z,w,\chi,&\tau,\bar
F(\chi,\tau),\bar G(\chi,\tau),(\partial
^\alpha\bar H(\chi,\tau))', (\partial
^\alpha\bar
H(\chi,\tau))''\bigr)\endaligned\tag3.4.26
$$
where $(\partial
^\alpha\bar H(\chi,\tau))=(\partial
^\alpha\bar H(\chi,\tau))_{|\alpha|\leq
\ell_0+|\nu|}$ and we use the
notation $(\partial
^\alpha\bar H(\chi,\tau))=((\partial
^\alpha\bar H(\chi,\tau))',(\partial
^\alpha\bar H(\chi,\tau))'')$ as above; 
\thetag{3.4.26} holds for
$w=Q(z,\chi,\tau)$ or
$\tau=\bar Q(\chi,z,w)$. 
We omit the details of this
construction, since it is similar to the
one above for the component $F$. Note, by
inspecting the construction above, that
the function $\Psi^{\tilde \jmath}_\nu$ is
defined only in terms of the defining
equations of
$M$ and
$M'$, and does not depend on the existence
or choice of a mapping $H$. The proof of
the formal part of the theorem is
complete.

Suppose that $M$ and $M'$ are also
real-analytic. Then the function
$S(\chi',\tau',r)$, defined by
\thetag{3.4.15}, is holomorphic in a
neighborhood of the origin. The fact,
noted above,  that each rational function
$R_{\alpha}
\left((a_\gamma), (b^{l}_\beta)\right)$
vanishes when $(b^{l}_\beta)=0$ implies
that the functions $\Psi^{\tilde
\jmath}_\alpha(Z,\zeta,\zeta',\Lambda)$ above
are holomorphic in a neighborhood of
$(Z,\zeta,\zeta';\Lambda',\Lambda'')=
(0,0,0,\Lambda_0',0)$ for any $\Lambda_0'$
such that \thetag{3.4.10} holds. This
completes the proof of Theorem 3.4.1.  
\qed\enddemo

We conclude this section by giving the
proof of Theorem 2.1.1.

\demo{Proof of Theorem $2.1.1$}  We take
normal coordinates for $M$ and $M'$ as in
the proof of Theorem 3.4.1. Let
$\tilde \jmath=(j_1,\ldots, j_{n'})$ be such
that $\det(\partial F^m_l/\partial
z_{j_p}(0,0))_{1\leq l,p\leq n'}\neq 0$, 
for $m=1,2$. By Proposition 3.1.6 and the
basic identity, Theorem 3.4.6, it follows
that
$$
\multline
\partial^\alpha
H(Z)\sim \Psi^{\tilde \jmath}_{\alpha}\left
(Z,\zeta,\bar H(\zeta),(\partial^\beta\bar
H(\zeta))_{1\leq |\beta|\leq
\ell_0+|\alpha|}\right),\\ \text{\rm for
} Z=v^{k+1}(z,\chi^1,z^1,\ldots),\
\zeta=\bar
v^k(\chi^1,z^1,\ldots),\endmultline
\tag3.4.27
$$ 
for any $k\geq 0$, where $v^l$ denotes the
Segre mapping defined in \S3.1 and
$v^0=(0,0)$. Hence, if
$j^{\ell_0k_0}_0(H^1)=j^{\ell_0k_0}(H^2)$,
then it follows from \thetag{3.4.27}, for
any
$k\leq k_0$, that
$$
(\partial^\alpha H^1)\circ v^k\sim
(\partial^\alpha H^2)\circ v^k,\ \forall
\alpha\: |\alpha|\leq
\ell_0(k_0-k),\tag3.4.28
$$
as can be seen by an induction on $k$.
In particular, we have 
$$
(H^1- H^2)\circ
v^{k_0}\sim 0.\tag3.4.29
$$
By Theorem 3.1.9, there exists $k_1$, with $k_1\leq d+1$, such that
$\Rk (v^{k})= N$, for $k\geq k_1$. It then follows from
standard facts about formal power series (see e.g.\ [BER4,
Proposition 5.3.5]) that
\thetag{3.4.29} implies $H^1\sim H^2$ if $k_0\geq k_1$. The proof is
complete.\qed 
\enddemo

\heading \S4. Proofs of the main results\endheading

\subhead \S4.1. Proof of Theorem
2.1.5\endsubhead It suffices to prove
Theorem 2.1.5 in 
normal coordinates. Thus, we take normal
coordinates for
$M$ and $M'$ as in previous sections. We
also keep the notation introduced in the
beginning of \S3.4 and in Remark
3.4.11. Consider the linear mapping
$D_k\:\bC^{kn}\to
\bC^{2kn}$ defined as follows. For $k=2j$,
$j\geq 1$, we set 
$$
\multline
D_{2j}(\chi^1,z^1,\ldots,
z^{j-1},\chi^j,z^j)=\\(0,\chi^1,z^1,\ldots,
z^{j-1},\chi^j,z^j,\chi^j,z^{j-1},\ldots,
z^1,\chi^1),\endmultline\tag4.1.1
$$
and for
$k=2j-1$, $j\geq 1$, we set
$$
D_{2j-1}(\chi^1,z^1,\ldots,
z^{j-1},\chi^j)=(0,\chi^1,z^1,\ldots,
z^{j-1},\chi^j,z^{j-1},\ldots,
z^1,\chi^1).\tag4.1.2
$$
We remind the reader that $v^k$ denotes
the $k$th Segre mapping as defined by
\thetag{3.1.4} and \thetag{3.1.5}. We
shall need the following.
\proclaim{Lemma 4.1.3} For any $k\geq 1$,
the following hold.
$$
v^{2k}\circ D_k\sim 0.\tag4.1.4
$$
For $k=2j$, $j\geq 1$,
$$
\rk \left(\frac{\partial
v^{2k}}{\partial z}\circ
D_k,\frac{\partial v^{2k}}{\partial
\chi^{j+1}}\circ D_k,\frac{\partial
v^{2k}}{\partial z^{j+1}}\circ
D_k,\ldots,\frac{\partial v^{2k}}{\partial
\chi^{k}}\circ D_k\right)=\Rk\,
(v^k).\tag4.1.5
$$
and 
for $k=2j-1$, $j\geq 1$,
$$
\multline
\rk \left(\frac{\partial
v^{2k}}{\partial z}\circ
D_k,\frac{\partial
v^{2k}}{\partial z^{j}}\circ
D_k,\frac{\partial v^{2k}}{\partial
\chi^{j+1}}\circ D_k,\frac{\partial
v^{2k}}{\partial z^{j+1}}\circ
D_k,\ldots,\frac{\partial v^{2k}}{\partial
\chi^{k}}\circ D_k\right)=\\\Rk\,
(v^k).\endmultline\tag4.1.6
$$
In particular, 
$$
\rk \left(\frac{\partial
v^{2k}}{\partial z}\circ
D_k,\frac{\partial v^{2k}}{\partial
\chi^{1}}\circ D_k,\frac{\partial
v^{2k}}{\partial z^{1}}\circ
D_k,\ldots,\frac{\partial v^{2k}}{\partial
\chi^{k}}\circ D_k\right)\geq\Rk\,
(v^k).\tag4.1.7
$$
\endproclaim 
\demo{Proof} Property \thetag{4.1.4}
follows by making repeated use of the
identities,
$$
Q(z,\chi,\bar Q(\chi,z,w))\sim w,\quad 
\bar Q(\chi,z,Q(z,\chi,\tau))\sim
\tau,\tag4.1.8
$$
which are easily
checked (see also [BER4, Remark 4.2.30]).

To prove \thetag{4.1.5} and \thetag{4.1.6},
we first write
$$
v^l(z,\chi^1,z^1\ldots,)=
(z,u^l(z,\chi^1,z^1\ldots)),\tag4.1.9
$$
where $u^l=(u^l_1,\ldots, u^l_d)$. We also
write $(z,\xi^{(l)})$ for the
variables
$(z,\chi^1,z^1\ldots)\in\bC^{ln}$. Observe,
by the form of
$v^l$ given by \thetag{3.1.4} and
\thetag{3.1.5}, that
$$
\rk\, \left(\frac{\partial
v^l(z,\xi^{(l)})}{\partial
(z,\xi^{(l)})}\right) = n+\rk\,\left(
\frac{\partial u^l(z,\xi^{(l)})}{\partial
\xi^{(l)}}\right).\tag4.1.10
$$
We shall complete the proof of
Lemma 4.1.3 in the case where $k=2j$,
and leave the odd case to the reader.
Thus, we shall prove \thetag{4.1.5}. We
have
$$
\multline
u^{4j}(z,\chi^1,z^1,\ldots,
\chi^{2j})=\\Q(z,\chi,\bar
Q(\chi,z^1,\ldots\bar
Q(\chi^j,z^j,u^{2j}(z^j,\chi^{j+1}\ldots,
z^{2j-1},\chi^{2j}))\ldots)).\endmultline
\tag 4.1.11
$$
For fixed $k=2j$, we shall write
$\xi^{(2k)} = (\xi',z^j,\xi'')$, where
$\xi'=(\chi^1,z^1,\ldots,\chi^j)$ and
$\xi''=(\chi^{j+1}, z^{j+1},\ldots,
\chi^{2j})$. We claim that 
$$
\rk\,\left(\frac{\partial
u^{2k}(z,\xi',z^j,\xi'')}{\partial\xi''}\circ
D_k
\right)=
\rk\,\left(\frac{\partial
u^{k}(z^j,\xi'')}{\partial\xi''}\right).
\tag4.1.12
$$ 
Since
$D_k(\xi',z^j)=(0,\xi',z^j,\tilde \xi')$,
where
$\tilde\xi'=(\chi^j,z^{j-1},
\ldots,\chi^1)$, \thetag{4.1.12} follows
from the chain rule, by using
\thetag{4.1.11} and the fact that
$$
\frac{\partial Q}{\partial
\tau}(0,0,0)=\frac{\partial \bar
Q}{\partial
w}(0,0,0)=I_{d\times d},\tag4.1.13
$$
where $I_{d\times d}$ denotes the
$d\times d$ identity matrix.
(The identity \thetag{4.1.13} is a
consequence of
\thetag{3.1.2}.) The desired
equality
\thetag{4.1.5} is an easy consequence of
\thetag{4.1.12}. This completes the proof
of Lemma 4.1.3.\qed\enddemo

We now return to the proof of Theorem
2.1.5.  We fix $\tilde \jmath$ as in that
theorem. Let $k_1$ be the integer provided by Theorem 3.1.9. We shall 
use the notation
$(z,\xi)$ for
$(z,\chi^1,z^1,\ldots)\in\bC^{2k_1n}$ as in
the proof of Lemma 4.1.3 above. We
claim that there exists a $\bC^{N'}$ valued
formal power series in $(z,\xi)$, 
$\Xi^{\tilde \jmath}(z,\xi,\Lambda')$, of the
form
$$
\Xi^{\tilde \jmath}(z,\xi,\Lambda')
\sim\sum_ {\gamma,\delta}\frac{
d_{\gamma\delta}(\Lambda')}
{\big(\det(\lambda^l_{j_p})_{1\leq
l,p\leq n'}\big)^ {l_{\gamma\delta}}}
z^\gamma\xi^\delta,
\tag4.1.14
$$
where
$d_{\gamma\delta}(\Lambda')$
are
$\bC^{N'}$ valued polynomials in
$\bC^{K'(2k_1\ell_0)}$
and
$l_{\gamma\delta}$ nonnegative integers,
satisfying the following.  For
every
$H\in\hat\Cal F(M,M')$ satisfying
\thetag{3.4.8}, we have
$$
H(v^{2k_1}(z,\xi))\sim\Xi^{\tilde
\jmath}\left (z,\xi,(\partial ^\alpha
H(0))'\right),\tag4.1.15
$$
where $(\partial ^\alpha
H(0))=(\partial ^\alpha
H(0))_{1\leq |\alpha|\leq
2k_1\ell_0}$ and $(\partial ^\alpha
H(0))=((\partial ^\alpha
H(0))',(\partial ^\alpha
H(0))'')$ as explained in Remark 3.4.11.
(Recall that $G(z,0)\sim 0$ so that
$\partial_z^\beta G(0,0)=0$ for all
$\beta$, i.e. $(\partial ^\alpha
H(0))''=0$.) Indeed, the existence of
such
$\Xi^{\tilde \jmath}$ follows by making
repeated use of
\thetag{3.4.27} for
$k=0,1,\ldots,2k_1-1$, complex conjugating every
other equation, and substituting
inductively. To check that $\Xi^{\tilde
\jmath}$ is of the form \thetag{4.1.14} is
straightforward, using \thetag{3.4.27}
and \thetag{3.4.7}, and is left to the
reader. 

In what follows, we shall assume that
$k_1=2j$ and leave the odd case to the
reader. As in the proof of Lemma 4.1.3, we
write $\xi=(\xi',z^j,\xi'')$, where
$\xi'=(\chi^1,z^1,\ldots,\chi^j)$ and 
$\xi''=(\chi^{j+1},z^{j+1},\ldots,
\chi^{2j})$. In view of Lemma 4.1.3, we
can choose $d$ components
$y''=(y''_1,\ldots, y''_d)$ from the
components of $\xi''$ such that
$$
\rk\,\left(\frac{\partial
u^{2k_1}(z,\xi',z^j,\xi'')}
{\partial y''}\circ D_{k_1}
\right)=d,\tag4.1.16
$$
where $u^{2k_1}$ is as defined in
\thetag{4.1.9}. After reordering the
components of $\xi''$ if necessary, we may
write $\xi''=(x'',y'')$, with
$x''=(x''_1,\ldots, x''_{(k_1-1)n-d})$
and $y''$ as above. We define the linear
isomorphism $m\:\bC^{(k_1-1)n}\to
\bC^{(k_1-1)n}$ so that
$$
D_{k_1}(\xi',z^j)=(\xi',z^j,m(\xi')).
\tag 4.1.17
$$
Using the decomposition
$\xi''=(x'',y'')$, the mapping
$m$ splits in the obvious way as
$m=(m_{x},m_{y})$.  We shall need the
following version of the implicit function
theorem with singularities.

\proclaim{Proposition 4.1.18} Let
$u(x,t,y)$ be a formal mapping 
$(\bC^{r_1}\times \bC^{r_2}\times
\bC^d,0)\to (\bC^d,0)$
such that 
$$
u(x,0,0)\sim 0,\quad
\rk\,\left(\frac{\partial
u}
{\partial y}(x,0,0)
\right)=d.\tag4.1.19
$$
Then the equation 
$$
u(x,t,y)\sim w,\tag4.1.20
$$
has a unique solution of the form
$$
y=\Delta(x)\theta\left(x,
\frac{t}{\Delta(x)^2},\frac{w}{\Delta(x)^2}
\right),\tag4.1.21
$$
where $\theta(t_1,t_2,t_3)$ is a formal
mapping $(\bC^{r_1}\times \bC^{r_2}\times
\bC^d,0)\to (\bC^d,0)$ and 
$$
\Delta (x)=\det \left(\frac{\partial
u}
{\partial y}(x,0,0)
\right).
\tag4.1.22
$$
If, in addition, the mapping $u$ is
holomorphic in a neighborhood of the
origin, then the mapping $\theta$ is also
holomorphic in a neighborhood of the
origin and \thetag{4.1.21} solves the
holomorphic equation $u(x,t,y)=w$ for
$(x,t,y,w)$ such that $\Delta(x)\neq 0$
and $|t/\Delta(x)^2|+|w/\Delta(x)^2|$ 
sufficiently small.
\endproclaim
\demo{Proof} It follows from the first
condition in \thetag{4.1.19} that
$$
u(x,t,y)=a(x,t,y)t+g(x,t,y)y,\tag4.1.23
$$
where $a(x,t,y)$ is a $d\times r_2$ matrix
of formal power series, $g(x,t,y)$ is a
$d\times d$ matrix of formal power series.
By expanding $g(x,t,y)$ in $t$ and $y$, we
obtain
$$
u(x,t,y)=g(x,0,0)y+y^\tau
R(x,t,y)y+\tilde a(x,t,y)t,\tag4.1.24
$$
where $R(x,t,y)$ is a $d\times d$
matrix of formal power series and
$\tilde a(x,t,y)$ is an $r_2\times d$
matrix of  formal power series. Note that
$$
g(x,0,0)\sim \frac{\partial
u}
{\partial y}(x,0,0).\tag4.1.25
$$ 
Using Cramer's rule on the equation
$u(x,t,y)\sim w$, we obtain, for some
$d\times d$ matrix $b(x)$ of formal power
series
$$
\Delta(x)y+b(x)y^\tau
R(x,t,y)y+b(x)\tilde a(x,t,y)t\sim b(x)
w,\tag4.1.26
$$
or, after dividing by $\Delta(x)^2$,
$$
\frac{y}{\Delta(x)}+b(x)\frac{y^\tau}
{\Delta(x)} 
R(x,t,y)\frac{y}{\Delta(x)}+b(x)\tilde
a(x,t,y)\frac{t}{\Delta(x)^2}\sim b(x)
\frac{w}{\Delta(x)^2}\tag4.1.27
$$
in $\bC[[x,t,y,w,1/\Delta(x)]]$. Let us
set $y'=y/\Delta(x)$, $t'=t/\Delta(x)^2$,
$w'=w/\Delta(x)^2$, and consider the
equation
$$
\multline
y'+b(x)(y')^\tau 
R(x,\Delta(x)^2t',\Delta(x)y')y'\\+b(x)\tilde
a(x,\Delta(x)^2t',\Delta(x)y')t'\sim b(x)
w'.\endmultline\tag4.1.27
$$ 
This has a unique formal mapping solution
$y'=\theta(x,t',w')$, with
$$
\theta\:(\bC^{r_1}\times \bC^{r_2}\times
\bC^d,0)\to (\bC^d,0),
$$ 
by the
formal implicit function theorem. The
conclusion of the proposition in the formal
case follows by substituting for $y'$,
$t'$, and $w'$. 

If $u$ is holomorphic in a neighborhood of
$0$, then the $\bC^d$ valued function
$\theta(x,t',w')$ is also holomorphic in a
neighborhood of
$0$, and it is straightforward to verify
the last conclusion of Proposition 4.1.18. 
\qed\enddemo

We return again to the proof of Theorem 2.1.5. We may apply
Proposition 4.1.18 to the equation
$$
u^{2k_1}(z,\xi',z^j,m_x(\xi'), y'')\sim
w,\tag4.1.28
$$
with $x=(\xi',z^j)$,
$t=z$, and $y=y''-m_y(\xi')$, since the
conditions in \thetag{4.1.19} are
satisfied by Lemma 4.1.3. We conclude
that the equation \thetag{4.1.28} has a
solution of the form
$$
y''=\phi\left
(\xi',z^j,\frac{z}{\Delta(\xi',z^j)^2},
\frac{w}{\Delta(\xi',z^j)^2}\right),
\tag4.1.29
$$ 
where $\phi$ is a formal mapping
$(\bC^{(k_1+1)n+d},0)\to (\bC^d,0)$ and 
$$
\Delta(\xi',z^j)=\det\left(\frac{\partial
u^{2k_1}}
{\partial
y''}\right)(\xi',z^j,m(\xi')).\tag4.1.30
$$
Thus, we have
$$
v^{2k_1}\left(z,\xi',z^j,m_x(\xi'),\phi\left
(\xi',z^j,\frac{z}{\Delta(\xi',z^j)^2},
\frac{w}{\Delta(\xi',z^j)^2}\right)\right)\sim
(z,w).\tag4.1.31
$$
Substituting this in \thetag{4.1.15}, we
obtain
$$
\multline
H(z,w)\sim\\\Xi^{\tilde
\jmath}\left (z,\xi',z^j,m_x(\xi'),\phi\left
(\xi',z^j,\frac{z}{\Delta(\xi',z^j)^2},
\frac{w}{\Delta(\xi',z^j)^2}\right),(\partial
^\alpha
H(0))'\right),\endmultline\tag4.1.32
$$
where $(\partial
^\alpha H(0))=(\partial
^\alpha H(0))_{1\leq |\alpha|\leq
2k_1\ell_0}$ and $(\partial
^\alpha H(0))=((\partial
^\alpha H(0))',(\partial
^\alpha H(0))'')$ as explained in Remark
4.1.11.  Expanding \thetag{4.1.32} as a
formal power series in
$z$ and $w$, we obtain, using \thetag{4.1.14},
$$
H(z,w)\sim\sum_{\alpha,\beta}\frac{
R_{\alpha\beta}(\xi',z^j,(\partial
^\alpha H(0))')}
{\Delta(\xi',z^j)^{2(|\alpha|+|\beta|)}}
z^\alpha w^\beta,\tag4.1.33
$$
where $R_{\alpha\beta}(\xi',z^j,\Lambda')$
is of the form
$$
R_{\alpha\beta}(\xi',z^j,\Lambda') 
\sim\sum_{\gamma,\delta}\frac{
R_{\alpha\beta\gamma\delta}(\Lambda')}
{\big(\det(\lambda^l_{z_{j_p}})_{1\leq
l,p\leq n'}\big)^
{l_{\alpha\beta\gamma\delta}}}
(\xi')^\gamma (z^j)^\delta,\tag4.1.34
$$
for polynomials
$R_{\alpha\beta\gamma\delta}(\Lambda')$ on
$\bC^{K'(2k_1\ell_0)}$ and nonnegative
integers
$l_{\alpha\beta\gamma\delta}$; we have
suppressed the dependence on $\tilde \jmath$
above to simplify the notation.

We next construct a formal mapping
$\hat\Phi(Z,\Lambda)=\hat\Phi^{\tilde
\jmath}(Z,\Lambda)$, with
$Z=(z,w)$ and
$\Lambda=(\Lambda',\Lambda'')$, as
follows.  First set
$$
\Gamma(z,w,\xi',z^j,\Lambda')=
\sum_{\alpha,\beta}\frac{
R_{\alpha\beta}(\xi',z^j,\Lambda')}
{\Delta(\xi',z^j)^{2(|\alpha|+|\beta|)}}
z^\alpha w^\beta.\tag4.1.35
$$ 
Since $\Delta(\xi',z^j)\not\sim 0$, we may
choose $({\xi'}^0,{z^j}^0)$ so that the
formal power series in $t$,
$$
\tilde
\Delta(t)=\Delta(t{\xi'}^0,t{z^j}^0)\not\sim
0.\tag4.1.36
$$
We shall write 
$$
\tilde \Gamma(z,w,t,\Lambda')=
\Gamma(z,w,t{\xi'}^0,t{z^j}^0,\Lambda')\sim
\sum_{\alpha,\beta}\frac{
\tilde R_{\alpha\beta}(t,\Lambda')}
{\tilde \Delta(t)^{2(|\alpha|+|\beta|)}}
z^\alpha w^\beta.\tag4.1.37
$$
where $\tilde R_{\alpha\beta}(t,\Lambda')=
R_{\alpha\beta}
(t{\xi'}^0,t{z^j}^0,\Lambda')$ are power
series in $t$ whose coefficients are
rational functions of $\Lambda$ of the
form appearing in \thetag{4.1.34}. By
division, we may write
$$
\tilde R_{\alpha\beta}(t,\Lambda')\sim
T_{\alpha\beta}(t,\Lambda')\tilde
\Delta(t)^{2(|\alpha|+|\beta|)}+
r_{\alpha\beta}(t,\Lambda'),\tag4.1.37
$$ 
where $T_{\alpha\beta}(t,\Lambda')$ is a
unique power series in $t$ whose
coefficients are finite linear
combinations of the coefficients of $\tilde
R_{\alpha\beta}(t,\Lambda')$, and
$r_{\alpha\beta}(t,\Lambda')$ is a
polynomial of degree at most
$2(|\alpha|+|\beta|)-1$ in $t$ whose
coefficients coincide with the
corresponding coefficient of $\tilde
R_{\alpha\beta}(t,\Lambda')$. We decompose
$\tilde \Gamma$ as follows
$$
\aligned
\tilde \Gamma(z,w, &t,\Lambda')\sim  
\sum_{\alpha,\beta}
T_{\alpha\beta}(t,\Lambda')z^\alpha
w^\beta + \sum_{\alpha,\beta}\frac{
r_{\alpha\beta}(t,\Lambda')}
{\tilde \Delta(t)^{2(|\alpha|+|\beta|)}}
z^\alpha w^\beta\\
\sim &\,\,
\tilde
\Phi(z,w,\Lambda')+\sum_{\alpha,\beta}
tS_{\alpha\beta}(t,\Lambda')z^\alpha
w^\beta + \sum_{\alpha,\beta}\frac{
r_{\alpha\beta}(t,\Lambda')} {\tilde
\Delta(t)^{2(|\alpha|+|\beta|)}} z^\alpha
w^\beta,
\endaligned\tag4.1.38
$$
where, on the second line, we have decomposed 
$T_{\alpha\beta}(t,\Lambda')\sim
T_{\alpha\beta}(0,\Lambda')+
tS_{\alpha\beta}(t,\Lambda')$ and where
$$
\tilde
\Phi(z,w,\Lambda')=\sum_{\alpha,\beta}
T_{\alpha\beta}(0,\Lambda')z^\alpha
w^\beta.\tag4.1.39
$$
We set
$\hat \Phi(z,w,\Lambda):=\tilde
\Phi(z,w,\Lambda')$ (so that
$\hat\Phi(z,w,\Lambda)$ is independent of
$\Lambda''$). We claim that 
$\hat\Phi(z,w,\Lambda)$ is of the form
\thetag{2.1.6}. Indeed, this is easy to
check from the above and is left to the
reader. We also claim that for each
$H\in\hat\Cal F(M,M')$ with $\det(\partial
F_l/\partial z_{j_p}(0,0))_{1\leq l,p\leq
n'}\neq 0$ we have
$$
H(z,w)\sim  \hat \Phi(z,w,(\partial ^\alpha
H(0)))\sim \tilde \Phi(z,w,(\partial 
^\alpha H(0))'),\tag4.1.40
$$
where $(\partial ^\alpha
H(0))=(\partial ^\alpha
H(0))_{1\leq |\alpha|\leq k_1\ell_0}$ and
$(\partial ^\alpha
H(0))=((\partial ^\alpha
H(0))',(\partial ^\alpha
H(0))'')$ as in Remark 3.4.11. In view of
\thetag{4.1.33}, \thetag{4.1.35},
\thetag{4.1.37}, and \thetag{4.1.38}, we
have
$$
\multline
H(z,w)-\hat\Phi(z,w,(\partial ^\alpha
H(0)))\sim\\ \sum_{\alpha,\beta}
\left(tS_{\alpha\beta}(t,(\partial ^\alpha
H(0))')+
\frac{
r_{\alpha\beta}(t,(\partial ^\alpha
H(0))')} {\tilde
\Delta(t)^{2(|\alpha|+|\beta|)}}\right)
z^\alpha w^\beta.
\endmultline\tag4.1.41
$$
Note that each coefficient on the right
hand side of \thetag{4.1.41} is a Laurent series in $t$
without constant term. Since the left hand
side is independent of $t$, we conclude that
each coefficient on the right hand side
must be zero and, hence, \thetag{4.1.40}
holds. 

Assume now that $M$ and
$M'$ are real-analytic. An inspection of
the proof above, using the fact that the
functions $\Psi^{\tilde
\jmath}_\alpha(Z,\zeta,\zeta',\Lambda)$ of
Theorem 3.4.6 are holomorphic near every
$(Z,\zeta,\zeta',\Lambda',\Lambda'')=
(0,0,0,\Lambda_0',0)$ for each
$\Lambda'_0$ satisfying \thetag{3.4.10},
and the fact that $\hat\Phi(Z,\Lambda)$
is independent of $\Lambda''$, shows that
$\hat\Phi(Z,\Lambda)$ is holomorphic near
every $(0,\Lambda_0)$ with
$\Lambda_0=(\Lambda_0',\Lambda_0'')$ such
that $\Lambda_0'$ satisfies
\thetag{3.4.10}.

Now, if we define the linear mapping
$\pi_{\tilde \jmath}$ in the conclusion of
Theorem 2.1.5 to be
$$
\pi_{\tilde \jmath}(\Lambda)=
(\lambda^l_{z_{j_p}})_{1\leq l,p\leq
n'},\tag 4.1.42
$$
and the polynomial $P$ to be 
$$
P\big((\lambda^l_{z_{j_p}})_{1\leq l,p\leq
n'}\big)=
\det\big((\lambda^l_{z_{j_p}})_{1\leq
l,p\leq n'}\big),\tag4.1.43
$$
then the formal mapping
$\hat\Phi(Z,\Lambda)$ above satisfies all
the conclusions of Theorem 2.1.5, except
that $\hat\Phi(Z,\Lambda)$ is a function
of $\Lambda\in
J^{2k_1\ell_0}(\bC^{N},\bC^{N'})_{(0,0)}$
and $2k_1$ need not be $\leq d+1$
(although $k_1\leq d+1$ as noted in \S1.3).
We shall address this point, and complete
the proof in the next section.

\subhead \S4.2. Conclusion of the proof
of Theorem 2.1.5\endsubhead To complete
the proof of the theorem, we shall need
the following proposition.

\proclaim{Proposition 4.2.1} Let
$F\:(\bC^l,0)\to (\bC^k,0)$ be a formal
mapping of the form
$$
F(x)=(F^0_1(x)+O(\kappa_1+1),\ldots,
F^0_k(x) + O(\kappa_k+1)),\tag4.2.2
$$
where $F^0_j(x)$ is a homogeneous
polynomial of degree $\kappa_j\geq 1$ for
$j=1,\ldots, k$. Assume that
$\Rk (F^0)=k$, where
$F^0=(F^0_1,\ldots, F^0_k)$. Then, for every
$\alpha\in\Bbb Z^k_+$, there exists a
linear form $\Cal
P_\alpha\:\bC^{\sigma(\|\alpha\|)}\to \bC$,
where
$$
\|\alpha\|=\sum_{j=1}^k\kappa_j\alpha_j
\tag4.2.3
$$
and 
$\sigma(\nu)$ denotes the number of
$\beta\in \Bbb Z_+^l$ with $|\beta|\leq \nu$,
such that the following holds. For every
$$
g(y)\sim\sum_{\alpha} g_\alpha
y^\alpha\in\bC[[y_1,\ldots,
y_k]],\quad h(x)\sim \sum_{\beta} h_\beta
x^\beta\in\bC[[x_1,\ldots,
x_l]]\tag4.2.4
$$
such that
$$
g(F(x))\sim h(x),\tag4.2.5
$$
we have
$$
g_\alpha=\Cal
P_\alpha\big((h_\beta)_
{|\beta|\leq\|\alpha\|}\big),\quad
\forall \alpha\in\Bbb Z_+^k.\tag4.2.6
$$
Moreover, the coefficients of the linear
form $\Cal P_\alpha$ are polynomial
functions of the coefficients of $F$.  
\endproclaim 
\demo{Proof} We decompose $g$ into
{\it weighted} homogeneous terms,
$$
g(y)\sim\sum_{\nu=0}^\infty
g^\nu(y),\tag4.2.7
$$
where $g^\nu(y)$ is a weighted
homogeneous polynomial of degree $\nu$ 
with respect to the weights
$(\kappa_1,\ldots, \kappa_k)$, i.e.
$$
g^\nu(y)\sim\sum_{\|\alpha\|=\nu}
g_\alpha y^\alpha.\tag4.2.8
$$
We
also decompose
$h(x)$ into standard homogeneous terms
$$
h(x)\sim \sum_{\nu=0}^\infty
h^\nu(x),\tag4.2.9
$$
where $h_\nu(x)$ is a
homogeneous polynomial of degree $\nu$.
Composing $g$ with $F$ and identifying
terms of degree $\nu$ in \thetag{4.2.5},
we obtain
$$
g^\nu(F^0(x))=h^\nu(x)+\ldots,\tag4.2.8
$$
where the dots signify terms involving
$g_\alpha$ for $\|\alpha\|<\nu$. Consider
the linear mapping $T_\nu\:\Cal
H^\kappa_\nu[y]\to \Cal H_\nu[x]$, where 
$\Cal H^\kappa_\nu[y]$ denotes the space
of weighted (with respect to
$\kappa=(\kappa_1,\ldots,\kappa_k)$)
homogeneous polynomials in
$y$ of degree $\nu$ and
$\Cal H_\nu[x]$ the space
of 
homogeneous polynomials in
$x$ of degree $\nu$, defined by
$$
T_\nu(g_\nu)=g_\nu\circ F^0.\tag4.2.9
$$
The fact that $\Rk(F^0(x))=k$
implies that $T_\nu$ is injective for
each $\nu=0,1,\ldots$. Hence, it has a
left inverse $L_\nu:
\Cal H_\nu[x]\to \Cal H^\kappa_\nu[y]$. It
follows that
$$
g^\nu=L_\nu(h^\nu+\ldots).\tag4.2.10
$$
Since, as mentioned above, the dots
involve only $g_\alpha$ with
$\|\alpha\|<\nu$, the proof is easily
completed by induction on $\nu$. \qed
\enddemo 

We now return again to the proof of Theorem
2.1.5. We shall keep 
$\tilde \jmath$ fixed as in the previous
section. Recall, from the proof of Theorem
3.1.9 in \S3.3, the notation
$v^k_0(z,\chi^1,\ldots)$ for the lowest
order homogeneous terms (in each
component) of the Segre mapping
$v^k(z,\chi^1,\ldots)$. Since $M$ is of
finite type at $0$, an inspection of the
proof of Theorem 3.1.9 shows that there
exists an integer
$k_1\leq d+1$ (also called
$k_1$ in the proof of Theorem 3.1.9) such
that $\Rk\, (v^k_0)=N$ for $k\geq k_1$. In
what follows, we shall assume that
$k_1$ is even (and leave the odd case to
the reader), and denote by
$\xi=(\chi^1,z^1,\ldots)\in\bC^{(k_1-1)n}$.
The same argument used to obtain
\thetag{4.1.15} shows that there exists a
$\bC^{N'}$ valued
formal power series in $(z,\xi)$, 
$\Xi^{\tilde \jmath}(z,\xi,\Lambda')$ of the
form
$$
\Xi^{\tilde \jmath}(z,\xi,\Lambda')
\sim\sum_ {\gamma,\delta}\frac{
d_{\gamma\delta}(\Lambda')}
{\big(\det(\lambda^l_{j_p})_{1\leq
l,p\leq n'}\big)^ {l_{\gamma\delta}}}
z^\gamma\xi^\delta,
\tag4.2.11
$$
where
$d_{\gamma\delta}(\Lambda')$
are
$\bC^{N'}$ valued polynomials in
$\bC^{K'(k_1\ell_0)}$
and
$l_{\gamma\delta}$ nonnegative integers,
satisfying the following.  For
every
$H\in\hat\Cal F(M,M')$ satisfying
\thetag{3.4.8}, we have
$$
H(v^{k_1}(z,\xi))\sim\Xi^{\tilde
\jmath}\left (z,\xi,(\partial ^\alpha
H(0))'\right),\tag4.2.12
$$
where $(\partial ^\alpha
H(0))=(\partial ^\alpha
H(0))_{1\leq |\alpha|\leq
k_1\ell_0}$ and $(\partial ^\alpha
H(0))=((\partial ^\alpha
H(0)',(\partial ^\alpha
H(0)'')$ as explained in Remark 3.4.11. We
shall apply Proposition 4.2.1 with
$x=(z,\xi)$, $y=Z=(z,w)$,
$F(x)=v^{k_1}(z,\xi)$,
$g(y)=H(Z)$, and $h(x)=\Xi^{\tilde
\jmath}(z,\xi,(\partial ^\alpha
H(0))')$. Since $F^0=v^{k_1}_0$ has rank
$N$ by definition of $k_1$, it follows
that the hypotheses on $F$ in Proposition
4.2.1 is satisfied with $\kappa_j=1$ for
$j=1,\ldots, n$, and $\kappa_{n+j}=\mu_j$
for $j=1,\ldots, d$, where the $\mu_j$
denote the H\"ormander numbers with
multiplicity as defined in \S3.2. We
conclude, using \thetag{4.2.11}, that for
every
$\beta\in\Bbb Z^N_+$, there is a linear
form $\Cal P_\beta$ such that
$$
\partial^\beta H(0)= \Cal
P_\beta\left(\left(\frac{
d_{\gamma\delta}((\partial ^\alpha
H(0)')}{\big(\det(\partial
F_l/\partial z_{j_p}(0))_{1\leq l,p\leq
n'}\big)^
{l_{\gamma\delta}}}\right)_{|\gamma|+|\delta|\leq
\|\beta\|}\right),\tag4.2.13
$$
where $(\partial ^\alpha
H(0)')$ is as in \thetag{4.2.12} (in
particular, $|\alpha|\leq k_1\ell_0$).
Now, let us write
$\Lambda'=(\Lambda'_{\alpha})$, where
$\Lambda'_{\alpha}$ stands for the
$\partial^\alpha H(0)$  part
of the jet $j^{2k_1\ell_0}_0(H)$ which
appears in 
$\Lambda'$ (see Remark 3.4.11).  By
substituting
$$
\Lambda'_\beta=\Cal
P_\beta\left(\left(\frac{
d_{\gamma\delta}(\Lambda'_\alpha)_{|\alpha|
\leq k_1\ell_0}}
{\big(\det(\lambda^l_{z_{j_p}})_{1\leq
l,p\leq n'}\big)^
{l_{\gamma\delta}}}\right)_{|\gamma|+|\delta|\leq
\|\beta\|}\right),\tag4.2.14
$$
for
$k_1\ell_0< |\beta|\leq 2k_1\ell_0$ in 
\thetag{4.1.39}, we obtain the desired
formal mapping 
$\Phi^{\tilde \jmath}(Z,\Lambda)$ in Theorem
2.1.5. In the case where $M$ and $M'$ are
real-analytic, we leave it to the reader
to check that substituting \thetag{4.2.14}
in \thetag{4.1.39} (which in this case
is holomorphic) yields a holomorphic
mapping as described in the theorem. This
completes the proof of Theorem
2.1.5.\qed

\subhead \S4.3. Proofs of Theorems 4, 2.1.9, 2.1.12, and 2.1.14\endsubhead 

\demo{Proof of Theorem $2.1.9$} Let $k_1$, $P$, $\pi_{\tilde \jmath}$, $\Phi^
{\tilde \jmath}$, for $\tilde \jmath=(j_1,\ldots,j_{n'})$ with $1\leq j_1
<\ldots< j_{n'}\leq n$, be given by Theorem 2.1.5. Then $\Lambda_0\in 
J^{k_1\ell_0}(\bC^N,\bC^{N'})_{(0,0)}$ is in the image of \thetag{2.1.10}
if and only if $P(\pi_{\tilde \jmath}\circ j^{k_1\ell_0,1}_0(\Lambda_0))
\neq 0$, for some $\tilde \jmath$, and
$$
\align
\Lambda_0 &=j_0^{k_1\ell_0}\left(\Phi^{\tilde \jmath}
(\cdot,\Lambda_0)\right),
\tag 4.3.1\\ 
\rho'\left(\Phi^{\tilde \jmath}(Z,\Lambda_0),\overline
{\Phi^{\tilde \jmath}}(\zeta,\bar \Lambda_0)\right) &\sim
a(Z,\zeta)\rho(Z,\zeta),
\tag 4.3.2
\endalign
$$
for some $d\times d$ matrix $a(Z,\zeta)$ of formal power 
series, where $\rho$ and $\rho'$ are defining power series
for $M$ and $M'$ respectively. In view of
\thetag{2.1.6}, the equation
\thetag{4.3.1} is a finite set of polynomial equations on
$\Lambda_0$. Similarly, by \thetag{2.1.6} and  elementary
linear algebra, \thetag{4.3.2} is an infinite set of 
polynomial equations on $\Lambda_0$ and $\bar \Lambda_0$.
Thus, the  solutions $\Lambda_0$  to \thetag{4.3.1} and
\thetag{4.3.2} form a (possibly empty) real-algebraic 
subvariety $A_{\tilde \jmath}$ of
$J^{k_1\ell_0}(\bC^N,\bC^{N'})_{(0,0)}$, and  the image of
\thetag{2.1.10} coincides with
$A\setminus B$, where
$$
A=\bigcup_{\tilde \jmath}A_{\tilde \jmath},\quad B=\bigcap_{\tilde \jmath} \{
P(\pi_{\tilde \jmath}\circ j^{k_1\ell_0,1}_0(\Lambda))=0\}.\tag4.3.3
$$
To see that the image is totally real at each regular point, we pick a
regular point $\Lambda_0\in A_{\tilde \jmath}\setminus B$ for some $\tilde
\jmath$. Thus, there is a unique mapping $H^0\in\hat\Cal F(M,M')$ with
$j^{k_1\ell_0}_0(H^0)=\Lambda_0$.
By applying the basic identity Theorem 3.4.6, with $Z=0$ and $\zeta=0$, to
the mapping $H^0$ and using \thetag{4.2.13} (cf. [BER3, Lemma
3.7]), we deduce that there exists a rational mapping $T(\Lambda)$,
holomorphic near $\Lambda_0$, such that
$$
\Lambda=T(\bar \Lambda),
\tag4.3.4
$$ 
holds for each $\Lambda\in A_{\tilde \jmath}\setminus B$ near $\Lambda_0$.
From this, it is easy to see that $A_{\tilde \jmath}\setminus B$ is totally
real at $\Lambda_0$. This completes the proof of Theorem 2.1.9.\qed
\enddemo

\demo{Proof of Theorem $4$} We shall prove Theorem 4 with $(d+1)\ell_0$ 
replaced 
by $k_1\ell_0$, where $k_1$ is as above. Clearly, this implies Theorem 4, 
since $k_1\leq d+1$. 
In view of the proof above and Theorem 3, it
suffices to show that the mapping \thetag{2.1.10} is a homeomorphisms onto
$A\setminus B$, with $A$ and $B$ as in the proof of Theorem 2.1.9, when $M$
and $M'$ are real-analytic. This is an easy consequence of Theorem 2.1.5
and the details are left to the reader (see also the proof of [BER3, 
Theorem 1]). \qed\enddemo

\demo{Proofs of Theorems $2.1.12$ and $2.1.14$} The conclusions of these theorems
follow immediately from Theorems 2.1.9 and 4 (or, more precisely, the version
of Theorem 4 with $(d+1)\ell_0$ replaced by $k_1\ell_0$ as proved above), 
since
a locally closed subgroup of a Lie group is a Lie subgroup (see e.g\
[Va]). \qed\enddemo

\heading \S 5. Smooth perturbations of formal generic
submanifolds\endheading

\subhead \S5.1. Smooth families of submanifolds\endsubhead
In this section, we shall study the behavior of the series
$\Phi^{\tilde \jmath}$, given in Theorem 2.1.5, under smooth perturbations of
the formal submanifolds $M$ and $M'$. For this, we need the notion of
a smooth family of formal generic 
submanifolds which will be introduced below. A particularly important
example, discussed in detail in \S5.2, is the family  
obtained by taking the 
formal generic
submanifold associated to a given smooth generic submanifold $M$ 
at varying points $p\in M$.

Let
$\rho(Z,\zeta;x)=(\rho_1(Z,\zeta;x),\ldots,\rho_d(Z,\zeta;x))$, 
where $Z=(Z_1, 
\ldots, Z_N)$ and $\zeta=(\zeta_1,\ldots, \zeta_N)$, be a 
smooth family of formal defining series, i.e. each 
$\rho_j$ is a formal power series in $Z$ and $\zeta$ whose coefficients are 
smooth functions of $x$ for $x$ in some smooth manifold $X$ and, for each 
fixed $x\in X$, $\rho(\cdot,\cdot;x)$ is a defining series of a formal
generic submanifold, denoted by $M_x$, as explained in \S1. 
The collection $\{M_x\}$, $x\in X$, will
be referred to as a {\it smooth family} of formal generic
submanifolds through
$0$ in $\bC^N$. If $X$ is a real-analytic manifold and the coefficients of 
$\rho$ depend real-analytically
on $x$, then we say that the family is real-analytic.
We have the 
following result. 

\proclaim{Theorem 5.1.1} Let $\{M_x\}$, $x\in X$, and $\{M'_y\}$, $y\in Y$, 
be smooth families of formal generic submanifolds of codimension $d$ and 
$d'$ through $0$ in $\bC^N$ and $\bC^{N'}$, respectively.  
Assume that $M_{x_0}$, for $x_0\in X$, is of finite type at $0$, that
$M'_{y_0}$, for $y_0\in Y$, is $\ell_0$-nondegenerate for some nonnegative
integer $\ell_0$, and that 
$n\geq n'$, where
$n=N-d$ and $n'=N'-d'$. 
Then there exist open neighborhoods $U\subset X$ and 
$V\subset Y$ of $x_0$ and $y_0$ respectively, 
an integer $k_1$ with
$1 < k_1\leq d+1$, polynomials
$P(\cdot;x,y)$ on $J^1(\bC^{n'},\bC^{n'})_{(0,0)}$ whose
coefficients depend smoothly  on 
$x,y\in U\times V$, and for
each $\tilde \jmath=(j_1,\ldots, j_{n'})$ with
$1\leq j_1<\ldots <j_{n'}\leq n$, linear
surjective mappings 
$$
\pi_{\tilde \jmath}(\cdot;x,y)\:
J^1(\bC^N,\bC^{N'})_{(0,0)}\to
J^1(\bC^{n'},\bC^{n'})_{(0,0)},\tag 5.1.2
$$
depending smoothly on $x,y\in U\times V$,
and a formal power series in
$Z=(Z_1,\ldots Z_N)$
$$
\Phi^{\tilde
\jmath}(Z,\Lambda;x,y)\sim\sum_{|\alpha|>0}\frac{
c^{\tilde \jmath}_\alpha(\Lambda;x,y)}{(P(\pi_{\tilde
\jmath}(j_0^{k_1\ell_0,1}(\Lambda);x,y);x,y))^
{l^{\tilde \jmath}_\alpha}} Z^\alpha, \tag5.1.3
$$
where $c^{\tilde \jmath}_\alpha(\cdot;x,y)$ are $\bC^{N'}$
valued polynomials on
$J^{k_1\ell_0}(\bC^N,\bC^{N'})_{(0,0)}$ whose coefficients depend smoothly 
on $x,y\in U\times V$ and
$l^{\tilde \jmath}_\alpha$ are nonnegative integers,
with the following property. For every formal CR submersive
mapping 
$H\:(M_x,0)\to(M_y',0)$, with $x,y\in U\times V$,
there exists $\tilde \jmath$ as above such that
$P(\pi_{\tilde \jmath}(j_0^{1}(H);x,y);x,y)\neq 0$ and 
$$
\aligned
H(Z)&\sim \Phi^{\tilde
\jmath}(Z,j^{k_1\ell_0}_0(H);x,y),\ \text{\rm if
$k_1$ is even,}\\ H(Z)&\sim \Phi^{\tilde
\jmath}\left(Z,\overline{j^{k_1\ell_0}_0(H)};x,y
\right),\ \text{\rm if
$k_1$ is odd}.
\endaligned\tag5.1.4
$$
In addition, if $M_x$ and $M'_y$, for some $x,y\in U\times V$, are
real-analytic, then for every $\tilde \jmath$
as above and $\Lambda_0$ with 
$P(\pi_{\tilde
\jmath}(j_0^{k_1\ell_0,1}(\Lambda_0);x,y);x,y)\neq 0$
the series \thetag{5.1.3} converges
uniformly for $(Z,\Lambda)$ near
$(0,\Lambda_0)$ in
$\bC^N\times
J^{k_1\ell_0}(\bC^N,\bC^{N'})_{(0,0)}$. If the families $\{M_x\}$, $x\in U$, 
and $\{M'_y\}$, $y\in Y$,
are real-analytic, then
the dependence of $P(\cdot;x,y)$, $\pi_{\tilde \jmath}(\cdot;x,y)$, and 
$c^{\tilde \jmath}_\alpha(\cdot,\cdot;x,y)$ above on $x,y\in U\times V$ is 
real-analytic for $x,y\in U\times V$.
\endproclaim

\demo{Proof} For the proof of Theorem 5.1.1, we shall need the
following version of the formal implicit function theorem in which the
dependence of the coefficients in the solution on the
coefficients in the equation is described. The proof, which consists of
applying the usual (formal) implicit function theorem and identifying
coefficients in the equation, is left to the reader.

\proclaim{Lemma 5.1.5} Let $k$ and $m$ be nonnegative integers. Then
there exist polynomials $\Cal
P_\gamma$, for $\gamma\in\Bbb Z_+^k$, with the following property.
For every formal power series mapping $F\:(\bC^k\times\bC^m,0) \to (\bC^m,0)$,
$$
F(x,y)\sim\sum_{\alpha,\beta}a_{\alpha\beta}x^\alpha y^\beta,\tag5.1.6
$$
with $x=(x_1,\ldots, x_k)$ and
$y=(y_1,\ldots, y_m)$, 
such that $\partial F/\partial
y(0,0)$ is invertible, there is a unique
power series solution $y=f(x)$ of the equation $F(x,y)\sim
0$, where $f\:(\bC^k,0)\to (\bC^m,0)$ is of the form
$$
f(x)\sim\sum_{\gamma}b_\gamma x^\gamma\tag5.1.7
$$
with 
$$
b_\gamma=\frac{\Cal
P_\gamma\left((a_{\alpha\beta})_{|\alpha|+|\beta|\leq|
\gamma|}\right)}{\det((a^j_k)_{1\leq j,k\leq
m})^{|\gamma|}},\quad\forall
\gamma\in \Bbb Z^k_+,
\tag5.1.8
$$
where $a^j_k=\partial F_j/\partial y_k(0,0)$. 
\endproclaim

We return to the proof of Theorem 5.1.1.
A consequence of Lemma 5.1.5 (and the construction of normal
coordinates; see [CM], [BJT], or [BER4, Chapter IV.2]) is that given a
smooth family $\{M_x\}$, 
$x\in X$, of formal generic submanifolds through $0$ in $\bC^N$ and
$x_0\in X$, there is a formal change to normal coordinates
$(z(Z;x),w(Z;x))$ whose coefficients depend smoothly on $x\in X$ near
$x_0$ such that $M_x$, in these coordinates, is defined by 
$w-Q(z,\chi,\tau;x)$, where $Q(z,\chi,\tau;x)$ is a $d$-vector of
formal power series in $(z,\chi,\tau)$ satisfying \thetag{1.12} and
whose coefficients depend smoothly on $x\in X$ near $x_0$. 
Now, Theorem 5.1.1 follows by a detailed inspection of the
proof of Theorem 2.1.5 and repeated use of Lemma 5.1.5. The
details are omitted. \qed 
\enddemo

It follows from Theorem 5.1.1 and the proof of Theorem
2.1.9 that, under the assumptions of Theorem
5.1.1, the defining equations of the images
$$
j^{k_1\ell_0}(\hat\Cal F(M_x,M'_y))\subset
J^{k_1\ell_0}(\bC^{N},\bC^{N'})
$$ 
depend smoothly on
$(x,y)\in X\times Y$ near $(x_0,y_0)$. Hence, by applying
[BER3, Lemma 5.1] (which essentially is the ``no small
subgroups'' property of Lie groups), we obtain in
particular the following corollary of Theorem 5.1.1.

\proclaim{Theorem 5.1.9} Let $\{M_x\}$, $x\in X$, be a smooth family
of formal generic submanifolds of codimension $d$ through $0$ in
$\bC^N$. Assume that $M_{x_0}$, for $x_0\in X$, is of finite type and 
$\ell_0$-nondegenerate, for some nonnegative
integer $\ell_0$, at $0$. Let $k_1$ be the integer obtained by
applying Theorem $2.1.12$ to $M_{x_0}$, and assume that
the Lie group $j^{k_1\ell_0}(\hat\Cal F(M_{x_0},M_{x_0}))\subset
G^{k_1\ell_0} (\bC^N)_0$ is discrete. Then there is an open
neighborhood $U\subset X$ of $x_0$ such that the Lie groups
$j^{k_1\ell_0}(\hat\Cal F(M_{x},M_{x}))\subset G^{k_1\ell_0}
(\bC^N)_0$, for $x\in U$, are also discrete.
\endproclaim

Theorem 5.1.9 implies e.g.\ that if $\hat \Cal
F(M_{x_0},M_{x_0})$ consists of only the identity, then for
$x\in X$ near $x_0$ there are no formal automorphisms ``near
the identity'' (in the sense of jets of order $k_1\ell_0$)
in $\hat \Cal F(M_{x},M_{x})$. If the family $\{M_x\}$
consists of real-analytic submanifolds, then Theorem 5.1.9
(combined with Theorems 3 and 2.1.14) yields the following
result. If $\Aut(M_{x_0},0)$ is discrete (in the natural
topology), then $\Aut(M_{x},0)$ is discrete for all $x\in X$
near $x_0$. 

\subhead \S5.2. Dependence on the base point in a smooth generic
submanifold\endsubhead We now turn to the particular example of a
smooth family obtained by taking the formal generic submanifolds
associated to points $p$ on a given smooth generic submanifold
$M\subset \bC^N$. More precisely, let $M\subset\bC^N$ be a smooth
generic submanifold 
with defining function $\rho(Z,\bar Z)$ near a distinguished point 
$p_0\in M$. We shall consider the smooth family $\{M_p\}$, for $p$ in a neighborhood 
$U$ of $p_0$ in $M$, of formal generic 
submanifolds through $0$ in $\bC^N$ defined by the Taylor
series at
$0$ of the smooth defining functions 
$$
\rho(Z,\zeta;p):=\rho(Z+p,\zeta+\bar p),\quad p\in
U.\tag5.2.1
$$
For smooth
generic submanifolds $M\subset \bC^N$, $M'\subset\bC^{N'}$ and points
$p\in M$, $p'\in M'$, we denote by $\hat\Cal F(M,p;M',p')$ the set of
formal mappings $H\:(\bC^N,p)\to (\bC^{N'},p')$ such
that $Z\mapsto H(Z+p)-p'$ maps $M_p$ into
$M'_{p'}$ and is CR submersive. Here, a formal mapping
$H\:(\bC^N,p)\to (\bC^{N'},p')$ is such that the components of
$H=(H_1,\ldots, H_{N'})$ are formal power 
series in $Z-p$ with constant term $H(p)$ equal to $p'$. We
shall need some more notation.

We denote by $E(\bC^N,\bC^{N'})_{(Z,Z')}$ the set of germs at $Z$ of
holomorphic mappings $(\bC^N,Z)\to (\bC^{N'},Z')$ and by
$\hat E(\bC^N,\bC^{N'})_{(Z,Z')}$ the set of formal mappings
$(\bC^N,Z)\to (\bC^{N'},Z')$. We denote by $E(\bC^N,\bC^{N'})$
the disjoint union
$$ 
E(\bC^N,\bC^{N'}) = \bigcup_{Z\in
\bC^N,Z'\in\bC^{N'}}E(\bC^N,\bC^{N'})_{(Z,Z')}, \tag 5.2.2
$$ 
and use similar notation for $\hat E(\bC^N,\bC^{N'})$. We define
$\hat\Cal F(M,M')_{(U,U')}\subset \hat E(\bC^N,\bC^{N'})$, for open
subsets $U\subset M$ and $U'\subset M'$, to be
$$ 
\hat\Cal F(M,M')_{(U,U')} = \bigcup_{p\in
U,p'\in U'}\hat\Cal F(M,p;M',p'),\tag 5.2.3
$$ 
and $\Cal F(M,M')_{(U,U')}$ similarly. We equip
$E(\bC^N,\bC^{N'})_{(Z,Z')}$ with the natural inductive limit
topology, and $E(\bC^N,\bC^{N'})$ with the topology it
inherits by the trivialization
$$
E(\bC^N,\bC^{N'}) \cong \bC^N\times\bC^{N'}\times
E(\bC^N,\bC^{N'})_{(0,0)}\tag5.2.4
$$
defined by taking a germ $H$ at $Z_0$ with $H(Z_0)=Z_0'$ to
$(Z_0,Z'_0,H_0)\in \bC^N\times\bC^{N'}\times
E(\bC^N,\bC^{N'})_{(0,0)}$ , where
$$
H_0(Z)=H(Z+Z_0)-Z'_0.\tag5.2.5
$$
For a positive integer $k$, we shall denote by $J^k(\bC^N,\bC^{N'})$
the complex manifold of
$k$-jets of holomorphic mappings $\bC^N\to \bC^{N'}$,
i.e.
$$ 
J^k(\bC^N,\bC^{N'}) = \bigcup_{Z\in
\bC^N,Z'\in\bC^{N'}}J^k(\bC^N,\bC^{N'})_{(Z,Z')} \tag5.2.6
$$ 
where $J^k(\bC^N,\bC^{N'})_{(Z,Z')}$ denotes the space of
$k$-jets at $Z$ of holomorphic mappings taking $Z$ to
$Z'$. We have a similar trivialization 
$$
J^k(\bC^N,\bC^{N'}) \cong \bC^N\times\bC^{N'}\times
J^k(\bC^N,\bC^{N'})_{(0,0)},\tag5.2.7
$$
and we have 
the jet mapping $j^k\:\hat E(\bC^N,\bC^{N'})\to
J^k(\bC^N,\bC^{N'})$. (See e.g. [GG] or [BER4, Chapter XII];
cf. also \S2 above.) We shall use the trivialization
\thetag{5.2.7} and refer to $(Z,Z',\Lambda)\in
\bC^N\times\bC^{N'}\times 
J^k(\bC^N,\bC^{N'})_{(0,0)}$ as coordinates for
$J^k(\bC^N,\bC^{N'})$.

We shall also use the notation $J^k(\bC^N,\bC^{N'})_{(U,U')}\subset J^k(\bC^N,\bC^{N'})$ for the
submanifold defined by 
$$ 
J^k(\bC^N,\bC^{N'})_{(U,U')} = \bigcup_{p\in U,p'\in
U'}J^k(\bC^N,\bC^{N'})_{(p,p')}, \tag5.2.8
$$ 
where $U$ and $U'$ are open subsets of $M$ and
$M'$ respectively. Observe that $j^k$ maps $\hat\Cal
F(M,M')_{(U,U')}$ into $J^k(\bC^N,\bC^{N'})_{(U,U')}$.

We have the following corollary of
Theorem 5.1.1, whose proof is similar to those of 
Theorems 4 and 2.1.9, and is left to the reader.

\proclaim{Theorem 5.2.9} Let $M$ and $M'$
be smooth generic submanifolds through
$p_0\in\bC^N$ and $p'_0\in\bC^{N'}$
respectively, such that
$M$ is of finite type at $p_0$ and $M'$ is
$\ell_0$-nondegenerate at $p'_0$ for some
integer $\ell_0$. Then there are open neighborhoods $U\subset M$ and
$U'\subset M'$ of $p_0$ and $p'_0$ respectively, an integer
$k_1$ depending only on $M$ with $1<k_1\leq d+1$ where $d$
denotes the codimension of $M$, a finite collection
$b_1(p,p',\Lambda),\ldots, b_l(p,p',\Lambda)$ of polynomials in
$\Lambda$ and a countable
collection $\{a_j(p,p',\Lambda,\bar\Lambda)\}_{j=1}^\infty$ of polynomials in
$(\Lambda, \bar\Lambda)$ whose coefficients are smooth functions of $(p,p')\in
U\times U'$ with the following property. The
image of the mapping
$$
j^{k_1\ell_0}\:\hat\Cal F(M,M')_{(U,U')}\to
J^{k_1\ell_0}(\bC^{N},\bC^{N'})_{(U,U')}
\tag5.2.10
$$
coincides with the locally closed set
$$
\{a_j(p,p',\Lambda,\bar\Lambda)=0,\
j=1,2,\ldots\}\setminus\{b_k(p,p',\Lambda)=0,\
k=1,2,\ldots l\}.\tag5.2.11
$$
If, in addition,  $M$ and $M'$ are real-analytic, then the coefficients of
$a_j(p,p',\Lambda,\bar\Lambda)$ and $b_k(p,p',\Lambda)$ depend real-analytically
on $(p,p')\in U\times U'$, the sets 
$$\aligned
&A=\{a_j(p,p',\Lambda,\bar\Lambda)=0,\
j=1,2,\ldots\},\\&
B=\{b_k(p,p',\Lambda)=0,\ k=1,2,\ldots l\}\endaligned\tag5.2.12
$$ 
are real-analytic subvarieties of $J^{k_1\ell_0}(\bC^{N},\bC^{N'})_{(U,U')}$,
and the mapping \thetag{5.2.10} (where
$\hat\Cal F(M,M')_{(U,U')}=\Cal F(M,M')_{(U,U')}$ by Theorem $3$) is a
homeomorphism onto its image $A\setminus B$.
\endproclaim

\heading \S6. Remarks on the algebraic mapping 
problem\endheading

We shall conclude this paper by showing how that arguments presented
here can be applied to the algebraic mapping problem; see
[BR], [BER1],  [M],  [Z3], [CMS] for recent work. 
This problem consists, loosely speaking, of finding conditions
which imply that any holomorphic mapping $H\:\bC^N\to\bC^{N'}$, 
defined near some point $p_0\in\bC^N$ and mapping a given
real-algebraic submanifold $M\subset\bC^N$ with $p_0\in M$ into
another $M'\subset\bC^{N'}$, is algebraic. Recall that a 
holomorphic mapping is called algebraic if all its components 
are algebraic functions, and a real submanifold is called real
algebraic if it is contained in a real algebraic subvariety of
the same dimension (see also [BER4, Chapter V]). 

Since this problem is not the main topic of
the present paper, we shall refer the reader to the papers
mentioned above for a more detailed discussion and 
history of the problem. We give here only the following
new result.

\proclaim{Theorem 6.1} Let $M\subset\bC^N$ be a real algebraic,
generic, and connected submanifold which is of finite type at some
point. Let $M'\subset\bC^{N'}$ be a real algebraic, generic
submanifold which is holomorphically nondegenerate and of
finite type at every point. Then any holomorphic mapping $H\:
\bC^N\to\bC^{N'}$, which is defined near some point $p_1\in M$, maps
$M$ into $M'$, and which is CR submersive at some point $p_2\in M$, 
is algebraic.
\endproclaim

\demo{Proof} Let $H$ be a holomorphic mapping as in the theorem
above. We denote by $\Omega\subset\bC^N$ a neighborhood of the point
$p_1$ in which $H$ is holomorphic. It follows from the assumptions
in the theorem that $M\cap\Omega$ is of finite type outside a proper 
real-analytic subvariety. Since the restriction of $H$ to $M\cap\Omega$ is
CR submersive outside a proper real-analytic subvariety, it
follows that there is a point $p_3\in M\cap\Omega$ at which $M$ is
of finite 
type and $H$ is CR submersive. Moreover, since $M'$ is of
finite type at $H(p_3)\in M'$, it follows from Corollary 1.27 that $H$
maps a neighborhood $U\subset M$ of $p_3$ onto a neighborhood
$U'\subset M'$ of 
$H(p_3)$. We deduce, using [BER1, Proposition 1.3.1] (see
also [BER4, Theorem 11.5.1]) and the fact that $M'$ is
holomorphically nondegenerate, that there is a point 
$p_4\in U$ at which $M$ is of finite type, $H$ is CR submersive,
and for which $M'$ is finitely nondegenerate at $H(p_4)$. Thus,
Theorem 2.1.5 applies with $p_0=p_4$. Now,
an inspection of the proof of Theorem 2.1.5 shows that the 
functions $\Phi^{\tilde\jmath}$ are algebraic when $M$ and $M'$
are real algebraic, and hence it follows from \thetag{2.1.8}
that $H$ is algebraic. 
\qed\enddemo

\remark{Remark $6.2$} The crucial result used in the proof
above, after applying Corollary 1.27, is the
following consequence of Theorem 2.1.5 (in the algebraic
setting): {\it If
$M\subset\bC^N$ 
is a real algebraic, generic submanifold which is of finite
type at
$p_0\in M$ and $M'\subset \bC^{N'}$ is a real algebraic,
generic submanifold which is finitely nondegenerate at
$p_0'\in M'$, then any holomorphic mapping
$H\:(\bC^N,p_0)\to (\bC^{N'},p_0')$ which maps $M$ into $M'$
and is CR submersive at $p_0$ is algebraic.} We should 
point out that this result also follows from 
[Z3, Theorem 1.6]. However, the approach in [Z3] differs
from the one in this paper.
\endremark

\enddocument
\end